\documentclass[10pt]{article}
\usepackage{amssymb}
\usepackage{amsmath}
\usepackage[dvips]{graphics}
\usepackage{epsfig}

\begin{document}
\title{Analytic Solutions of the Heat Equation}

\author{\textbf{Vassilis G. Papanicolaou$^1$, Eva Kallitsi$^2$, and George Smyrlis$^3$}
\\\\
Department of Mathematics
\\
National Technical University of Athens,
\\
Zografou Campus, 157 80, Athens, GREECE
\\
{\tt $^1$papanico@math.ntua.gr, $^2$evpapa@hotmail.com,}
\\
{\tt $^3$gsmyrlis@math.ntua.gr}}
\maketitle

\begin{abstract}
Motivated by the recent proof of Newman's conjecture \cite{R-T} we study certain properties of entire caloric functions, namely solutions of the heat equation $\partial_t F = \partial_z^2 F$ which are entire in $z$ and $t$. As a prerequisite, we establish some general properties of the order and type of an entire function. Then, we start our inquiry on entire caloric functions by determining the necessary and sufficient condition for a function $f(z)$ to be the initial condition of an entire solutions of the heat equation and, subsequently, we examine the relation of the $z$-order and $z$-type of an entire caloric function $F(t, z)$, viewed as function of $z$, to its $t$-order and $t$-type respectively, if it is viewed as function of $t$. After that, we shift our attention to the zeros $z_k(t)$ of an entire caloric function $F(t, z)$, viewed as function of $z$. We
show that the points $(t, z)$ at which $F(t, z) = \partial_z F(t, z) = 0$ form a discrete set in $\mathbb{C}^2$ and we derive the
$t$-evolution equations of the zeros of $F(t, z)$. These are differential equations which hold for all but countably many $t \in \mathbb{C}$.
\end{abstract}

\textbf{Keywords.} Entire solutions of the heat equation; entire caloric functons; order; type; dynamics of the zeros.
\\\\
\textbf{2010 AMS Mathematics Classification.} 32A15; 32W30.

\section{Introduction}
In this paper we study certain properties of the solutions $F(t, z)$ of the standard heat equation $u_t = u_{zz}$
which are analytic
in $z$ and $t$. The motivation for our study sprang from the recent article \cite{R-T} by Rodgers and Tao on the proof of the Newman's conjecture.
This conjecture, which can be considered a complement to the Riemann hypothesis, appeared in Newman's paper \cite{N} (in Remark 2, p. 247).

Let us briefly summarize the main facts and ideas relevant to Newman's conjecture, in order to show the connection with the analytic solutions of the
heat equation. Following the notation of \cite{R-T} we set
\begin{equation}
H_0(z) := \frac{1}{8} \xi\left(\frac{1}{2} + \frac{iz}{2}\right),
\qquad\text{where }\;
\xi(s) := \frac{s(s-1)}{2} \pi^{-s/2} \Gamma\left(\frac{s}{2}\right) \zeta(s)
\label{RT1}
\end{equation}
($\zeta$ is the Riemann zeta function---$\xi$ is usually called the ``$\,$Riemann xi function"). The function $H_0$ is entire of order $1$ and of maximal (i.e. infinite) type, and satisfies $H_0(-z) = H_0(z)$, i.e. it is even (the evenness of $H_0$ is equivalent to the functional equation of the zeta function). The Riemann hypothesis, as it was originally stated by Riemann \cite{R}, is the conjecture that all the zeros of $H_0$ are real.
Moreover, Riemann \cite{R} had essentially derived the Fourier representation
\begin{equation}
H_0(z) = \int_0^{\infty} \Phi(x) \cos(zx) \, dx,
\label{RT2}
\end{equation}
where $\Phi$ is the super-exponentially decaying function (as $x \to \infty$)
\begin{equation}
\Phi(x) := \sum_{n=1}^{\infty} \left(2\pi^2 n^4 e^{9x} -3\pi n^2 e^{5x}\right) \exp\left(-\pi n^2 e^{4x}\right).
\label{RT3}
\end{equation}
The series in \eqref{RT3} converges absolutely for every real $x$ (and even for complex $x$ as long as $\cos(4 \, \Im(x)) > 0$) and, with the help of Poisson summation formula, it can be verified that $\Phi$ is, actually, even.

In a fruitful attempt to approach the Riemann hypothesis, de Bruijn \cite{dB} introduced the more general function
\begin{equation}
H(t, z) := \int_0^{\infty} e^{tx^2} \Phi(x) \cos(zx) \, dx.
\label{RT4}
\end{equation}
Obviously, $H(t, z)$ is even in $z$ and it is not really hard to show that it is entire in $z$ of order $1$ and of maximal type
(notice that it is entire in $t$ too). Furthermore, it is easy to check that $H$ satisfies the backwards heat equation
$\partial_t H = -\partial_z^2 H$ (thus $\hat{H}(t, z) := H(-t, z)$ satisfies the
standard heat equation), with $H(0, z) = H_0(z)$. A crucial property of $H$ is \cite{N} that there is a constant $\Lambda$, now known as the
\emph{de Bruijn-Newman constant}, with the property that $H(t, z)$ has purely real zeros if and only if
$t \geq \Lambda$ (here $t$ is viewed as a real parameter).
Thus, the Riemann hypothesis is equivalent to the upper bound $\Lambda \leq 0$. Newman \cite{N} conjectured the complementary lower bound $\Lambda \geq 0$, and noted that his conjecture asserts that if the Riemann hypothesis is true, it is only ``$\,$barely so."

After some considerable progress achieved in the earlier works \cite{C-S-V}, \cite{K-K}, \cite{K-K-L} (based on the study of the zeros of certain entire solutions of the heat
equation), Rodgers and Tao \cite{R-T} managed to prove the Newman's conjecture. One key ingredient in their proof is the system
of ordinary differential (evolution) equations satisfied by the zeros $z_k(t)$ of $H(t, z)$, namely
\begin{equation}
z'_k(t) = 2 \sum_{j \ne k} \frac{1}{z_k(t) - z_j(t)},
\qquad \text{where }\;
t > \Lambda.
\label{RT5}
\end{equation}
The above equations were first derived by Csordas, Smith, and Varga \cite{C-S-V}. They look like a kind of
``\,characteristics" for the heat equation and they also remind the equations which arise in the solution of the inverse spectral problem for the Hill operator (see, e.g., \cite{T}).

Thus, the theory developed for the establishment of the Newman's conjecture exhibited some interesting phenomena regarding analytic
solutions of the heat equation, and this suggested to us that the topic deserves an independent study.

\subsection{Entire solutions of the heat equation}
Suppose that the function $F(t, z)$ is entire in $z$ for every $t \in \mathbb{C}$ and entire in $t$ for every $z \in \mathbb{C}$. Then, a consequence
of the celebrated theorem of Hartogs (see, e.g., \cite{K}) is that, for any given pair $(t_0, z_0)$ of complex numbers the function $F(t, z)$ equals
to its Taylor expansion about $(t_0, z_0)$, namely
\begin{equation}
F(t, z) = \sum_{j,k \geq 0} c_{jk} (t - t_0)^j (z - z_0)^k,
\qquad \text{where} \quad
c_{jk} = \frac{\partial_t^j \partial_z^k F(t_0, z_0)}{j! k!}
\label{P1}
\end{equation}
and the series converges absolutely for any $t, z \in \mathbb{C}$. A further consequence of the absolute convergence is that $F(t, z)$ can be expanded as
\begin{equation}
F(t, z) = \sum_{k \geq 0} a_k(t) (z - z_0)^k,
\quad \text{as well as} \quad
F(t, z) = \sum_{j \geq 0} b_j(z) (t - t_0)^j,
\label{P2}
\end{equation}
where $a_k(t)$, $k \geq 0$, and $b_j(z)$, $j \geq 0$, are entire functions.

Let us now assume that $F(t, z)$ is also \emph{caloric}, namely it satisfies the heat equation
\begin{equation}
\partial_t F(t, z) = \partial_z^2 F(t, z)
\label{P3}
\end{equation}
with ``$\,$initial condition"
\begin{equation}
f(z) := F(0, z)
\label{P3a}
\end{equation}
(notice that if \eqref{P3} is satisfied in an open subset of $\mathbb{C}^2$, then, by analytic continuation of
$\partial_t F(t, z)$ and $\partial_z^2 F(t, z)$ we have that \eqref{P3} is automatically satisfied
for every $(t, z) \in \mathbb{C}^2$). Since the operators $\partial_t$ and $\partial_z$ commute, by differentiating \eqref{P3} with respect to $t$ repeatedly we obtain
\begin{equation}
\partial_t^j F(t, z) = \partial_z^{2j} F(t, z),
\qquad
j \geq 0.
\label{P4}
\end{equation}
Also, a rather trivial observation is that if $F(t, z)$ satisfies \eqref{P4}, so does
\begin{equation*}
\tilde{F}(t, z) = F(t + t', z + z').
\end{equation*}

Using \eqref{P4} in \eqref{P1} yields
\begin{equation}
F(t_0 + t, z_0 + z) = \sum_{j,k \geq 0} \frac{\partial_z^{2j+k} F(t_0, z_0)}{j! k!} \, t^j z^k,
\qquad
t, z \in \mathbb{C}.
\label{P5}
\end{equation}
It is more convenient to write \eqref{P5} as
\begin{equation}
F(t_0 + t, z_0 + z) = \sum_{m=0}^{\infty} \frac{\partial_z^m F(t_0, z_0)}{m!} \, P_m(t, z)
\label{P6}
\end{equation}
where
\begin{equation}
P_m(t, z) := \sum_{2j+k = m} \frac{m!}{j! k!} \, t^j z^k
= \sum_{j = 0}^{\lfloor m/2\rfloor} \frac{m!}{j! (m-2j)!} \, t^j z^{m-2j},
\qquad
m \geq 0.
\label{P7}
\end{equation}
The quantity $P_m(t, z)$ is called the \emph{$m$-th caloric polynomial} and it is clear from \eqref{P7} that it is
\emph{parabolically $m$-homogeneous}, namely
\begin{equation}
P_m(\lambda^2 t, \lambda z) = \lambda^m P_m(t, z).
\label{P8}
\end{equation}
The first six caloric polynomials are
\begin{equation*}
P_0(t, z) \equiv 1,
\qquad
P_1(t, z) = z,
\qquad
P_2(t, z) = z^2 + 2t,
\qquad
P_3(t, z) = z(z^2 + 6t),
\end{equation*}
\begin{equation*}
P_4(t, z) = z^4 + 12tz^2 + 12t^2,
\qquad
P_5(t, z) = z(z^4 + 20tz^2 + 60t^2).
\end{equation*}

Let us review some other important properties of the caloric polynomials. For each $m \geq 0$ the polynomial $P_m(t, z)$ satisfies the heat equation
with initial condition
\begin{equation}
P_m(0, z) = z^m,
\label{P8a}
\end{equation}
hence, the entire solution $F(t, z)$ of \eqref{P3}, with $f(z) = F(0, z) = \sum_{m=0}^M a_m z^m$,
is the polynomial $F(t, z) = \sum_{m=0}^M a_m P_m(t, z)$.

Also, from the standard integral formula, involving the heat kernel, which gives the solution of the heat equation in terms of the initial condition,
we have
\begin{equation}
P_m(t, z) = \int_{-\infty}^{\infty}\frac{1}{2\sqrt{\pi t}} \, e^{-(z - \xi)^2/4t} \, \xi^m d\xi
\qquad \text{for} \quad \Re(t) > 0.
\label{P8b}
\end{equation}
Since
\begin{equation*}
\frac{m!}{j! (m-2j)!} = \frac{m!}{(2j)! (m-2j)!} \, (j+1)(j+2) \cdots 2j = \binom{m}{2j} (j+1)(j+2) \cdots 2j
\end{equation*}
it follows that the coefficients of the caloric polynomials are positive integers.

It is easy to see from \eqref{P7} that a crude bound of
$P_m(t, z)$ is
\begin{equation}
|P_m(t, z)| \leq \frac{m!(\lfloor m/2\rfloor + 1)}{(\kappa_m)! (m-2\kappa_m)!} \,\,
\max_{0 \leq j \leq \lfloor m/2\rfloor} |t|^j |z|^{m-2j},
\label{P7b}
\end{equation}
where
\begin{equation}
\kappa_m = \left\lfloor \frac{4m-1-\sqrt{8m+17}}{8}\right\rfloor + 1.
\label{P7c}
\end{equation}

If we differentiate \eqref{P7} with respect to $z$ we get
\begin{equation}
\partial_{z} P_m(t, z) = \sum_{2j+k = m} \frac{m!}{j! (k-1)!} \, t^j z^{k-1}
= m \sum_{2j+l = m-1} \frac{(m-1)!}{j!\, l!} \, t^j z^l = m P_{m-1}(t, z)
\label{P7d}
\end{equation}
for $m \geq 1$ (another way to see that $\partial_{z} P_m(t, z) = m P_{m-1}(t, z)$ is by observing that, since $P_m(t, z)$ is the solution of the
heat equation \eqref{P3} with $P_m(0, z) = z^m$, the derivative $\partial_{z} P_m(t, z)$ is the solution of \eqref{P3} with initial condition
$m z^{m-1}$).

The function
\begin{equation}
E_{\lambda}(t, z) := e^{\lambda^2 t + \lambda z}
\label{P9}
\end{equation}
is entire in $(t, z)$ and satisfies the heat equation \eqref{P3}. Thus, we can apply \eqref{P6} to $E_{\lambda}(t, z)$ (for $t_0 = z_0 = 0$) and obtain
\begin{equation}
E_{\lambda}(t, z) = e^{\lambda^2 t + \lambda z} = \sum_{m=0}^{\infty} \frac{\partial_z^m E_{\lambda}(0, 0)}{m!} \, P_m(t, z)
= \sum_{m=0}^{\infty} \frac{\lambda^m}{m!} \, P_m(t, z).
\label{P10}
\end{equation}
In other words, $e^{\lambda^2 t + \lambda z}$ is the generating function of the caloric polynomials.

Formula \eqref{P7} also implies
\begin{equation}
P_m(-1, 2z) = m! \sum_{j = 0}^{\lfloor m/2\rfloor} \frac{(-1)^j}{j! (m-2j)!} \, (2z)^{m-2j} = H_m(z),
\qquad
m \geq 0,
\label{P11}
\end{equation}
where $H_m(z)$ is the (physicists') Hermite polynomial of order $m$, i.e. of degree $m$ (see, e.g., \cite{S}). Thus, the parabolic
homogeneity \eqref{P8} yields
\begin{equation}
P_m(t, z) = \left(i \sqrt{t}\right)^m H_m\left(\frac{z}{2i\sqrt{t}}\right)
\label{P11a}
\end{equation}
and equation \eqref{P6} can be written as (for $t_0 = z_0 = 0$)
\begin{equation}
F(t, z) = \sum_{m=0}^{\infty} \frac{\partial_z^m F(0, 0)}{m!} \, \left(i \sqrt{t}\right)^m H_m\left(\frac{z}{2i\sqrt{t}}\right)
= \sum_{m=0}^{\infty} \frac{f^{(m)}(0)}{m!} \, \left(i \sqrt{t}\right)^m H_m\left(\frac{z}{2i\sqrt{t}}\right),
\label{P6a}
\end{equation}
where the second equality follows from \eqref{P3a}. Thus, if $F(t, z)$ and $G(t, z)$ are two entire solutions of \eqref{P3} with
$F(0, z) \equiv G(0, z)$ (or $F(t_0, z) \equiv G(t_0, z)$ for some fixed $t_0 \in \mathbb{C}$), then \eqref{P6a} tells us that they have to be
identical, namely $F(t, z) \equiv G(t, z)$. Equivalently, if for some fixed $t_0$ we have $F(t_0, z) \equiv 0$, then $F(t, z) \equiv 0$. On the other hand, the relation $F(t, z_0) \equiv 0$, for some $z_0$ does not imply that $F$ is identically $0$ (e.g.,
if $f(z) \not\equiv 0$ is odd, then $F(t, z)$ is a solution of
\eqref{P3} which is odd in $z$, and hence $F(t, 0) \equiv 0$).

There are many known facts about the zeros of the Hermite polynomials. For instance, it is well known \cite{S} that $H_m(z)$ is an even (odd) function if and only if $m$ is even (odd). Furthermore, the zeros of $H_m(z)$ are real and simple. It, then, follows from \eqref{P11a} that \cite{P-S} if $m = 2l$, the polynomial $P_m$ is of the form
\begin{equation}
P_m(t, z) = (z^2 + \rho_{m,1} \, t) \cdots (z^2 + \rho_{m,l} \, t),
\qquad \text{with} \quad
0 < \rho_{m,1} < \cdots < \rho_{m,l},
\label{P12}
\end{equation}
while if $m = 2l + 1$, then $P_m$ is of the form
\begin{equation}
P_m(t, z) = z (z^2 + \rho_{m,1} \, t) \cdots (z^2 + \rho_{m,l} \, t),
\qquad \text{with} \quad
0 < \rho_{m,1} < \cdots < \rho_{m,l}.
\label{P13}
\end{equation}
From \eqref{P12} and \eqref{P13} we have that if $t \in \mathbb{C} \setminus \{0\}$, then the zeros of $P_m(t, z)$, viewed as a polynomial of $z$, are simple (the case $t=0$ is exceptional since $P_m(0, z) = z^m$). Furthermore, by
\eqref{P7d} and Rolle's Theorem we get the interlacing properties
\begin{equation}
0  < \rho_{m,1} < \rho_{m-1,1} < \rho_{m,2} < \cdots < \rho_{m-1,l-1} < \rho_{m,l},
\qquad \text{if} \quad
m = 2l,
\label{P12a}
\end{equation}
while
\begin{equation}
0 < \rho_{m-1,1} < \rho_{m,1} < \rho_{m-1,2} < \cdots < \rho_{m-1,l} < \rho_{m,l},
\qquad \text{if} \quad
m = 2l+1.
\label{P13a}
\end{equation}
Let us also notice that \eqref{P12} and \eqref{P13} tell us that, if $m > 0$ is even, then the zeros of $P_m(t, z)$ (viewed as a function of $z$) are
real if and only if $t \in \mathbb{R}^- := (-\infty, 0]$ and the same is true for the non-zero zeros of $P_m(t, z)$ in the case where $m$ is odd
($z = 0$ is always a zero of $P_m(t, z)$, if $m$ is odd).

The rest of the paper is organized as follows. In Section 2 we present some general results regarding the order and the type of an entire function. These results, mainly Theorems 1 and 2, apart from having their own interest, they will be used in Section 3, where we study the relation of the orders $\rho_z$ and $\rho_t$, as well as the types $\tau_z$ and $\tau_t$, of an entire caloric function $F(t, z)$, viewed as function of $z$ and $t$ respectively. The main results of Section 3 are Theorems 3, 4, and 5.
Finally, in Section 4 we first show that the multiple zeros of $F(t, z)$ are isolated (Theorem 6) and, then, we derive the dynamics of the zeros of
$F(t, z)$ for $(t, z) \in \mathbb{C}^2$.

\section{A general discussion on the order and the type of an entire function}
Let
\begin{equation}
g(z) = \sum_{n \geq 0} a_n z^n,
\qquad
z \in \mathbb{C},
\label{PREL1}
\end{equation}
be an entire function and
\begin{equation}
M(r) = M_g(r) := \sup_{|z| \leq r} |g(z)| = \max_{|z| = r} |g(z)|,
\qquad
r > 0,
\label{PREL2}
\end{equation}
its maximum modulus.

We recall that the order of $g(z)$ is the quantity \cite{H}
\begin{equation}
\rho = \rho(g) := \limsup_{r \to \infty} \frac{\ln\ln M(r)}{\ln r}.
\label{PREL3}
\end{equation}
In other words, the order $\rho$ of $g(z)$ is the smallest exponent $\rho' \geq 0$ such that for any given $\varepsilon > 0$ there is a
$r_0 = r_0(\varepsilon) > 0$ for which
\begin{equation}
|g(z)| \leq \exp\left(|z|^{\rho' + \varepsilon}\right)
\qquad \text{whenever }\;
r = |z| \geq r_0.
\label{PREL4}
\end{equation}
Clearly, $0 \leq \rho \leq \infty$.

Let us also recall \cite{H} that if $0 < \rho < \infty$, the quantity
\begin{equation}
\tau = \tau(g) := \limsup_{r \to \infty} \frac{\ln M(r)}{r^{\rho}}
\label{PREL5}
\end{equation}
is the type of (the order of) $g(z)$. In other words, $\tau$ is the smallest (extended) number $\tau' \geq 0$ such that for any given
$\varepsilon > 0$ there is a $r_0 = r_0(\varepsilon) > 0$ for which
\begin{equation}
|g(z)| \leq \exp\Big((\tau' + \varepsilon) |z|^{\rho}\Big)
\qquad \text{whenever }\;
r = |z| \geq r_0.
\label{PREL6}
\end{equation}
Clearly, $0 \leq \tau \leq \infty$. If $\tau = 0$, we say that $g(z)$ is of minimal type, whereas if $\tau = \infty$, we say that $g(z)$ is of
maximal type. In the extreme cases where $\rho = 0$ or $\rho = \infty$ the type is not defined.

It is clear from \eqref{PREL4} and \eqref{PREL6} that the entire functions
\begin{equation}
g(z)
\quad \text{and} \quad
\tilde{g}(z) := g(z + z_0)
\qquad\quad \text{have the same order and type}
\label{PREL7}
\end{equation}
for any fixed $z_0 \in \mathbb{C}$.

A well-known fact of complex analysis is \cite{H} that the order $\rho$ and the type $\tau$ of $g(z)$ are given by the formulas
\begin{equation}
\rho = \limsup_n \frac{n \ln n}{-\ln |a_n|}
\label{PREL8}
\end{equation}
and (in the case where $0 < \rho < \infty$)
\begin{equation}
\tau = \frac{1}{e \rho}\limsup_n n |a_n|^{\rho/n}
\label{PREL9}
\end{equation}
respectively, where $a_n$, $n = 0, 1, \ldots$, are the coefficients of the power series of $g(z)$ as seen in \eqref{PREL1}.

%
%

Another tool that we will need in the sequel is the operator $(^{\sharp})$ defined as
\begin{equation}
g^{\sharp}(z) := \sum_{k \geq 0} |a_k| z^k,
\label{PREL9a}
\end{equation}
when $g(z)$ is the entire function of \eqref{PREL1}. Notice that for any $r > 0$ we have
\begin{equation}
\max_{|z| \leq r} \left|g(z)\right| \leq g^{\sharp}(r) = \max_{|z| \leq r} \left|g^{\sharp}(z)\right|
\label{PREL9b}
\end{equation}
(the inequality can be strict).
Furthermore, since the order and type of $g(z)$ depend only on $\{|a_k|\}_{k \geq 0}$, they remain invariant under $(^{\sharp})$, i.e.
\begin{equation}
\rho\left(g^{\sharp}\right) = \rho(g)
\qquad \text{and} \qquad
\tau\left(g^{\sharp}\right) = \tau(g).
\label{PREL9c}
\end{equation}
Also,
\begin{equation}
\left(g'\right)^{\sharp}(z) = \sum_{k \geq 0} |k a_k| z^{k-1} = \sum_{k \geq 0} k |a_k| z^{k-1} = \left(g^{\sharp}\right)'(z),
\label{PREL9d}
\end{equation}
i.e. $(^{\sharp})$ commutes with the derivative operator.

Let us now set
\begin{equation}
a_n(z) := \frac{g^{(n)}(z)}{n!}
\qquad
n = 0, 1, \ldots
\label{PREL10}
\end{equation}
(so that $a_n(0) = a_n$). Then, in view of \eqref{PREL7}, formulas \eqref{PREL8}, \eqref{PREL9}, and \eqref{PREL10} yield
\begin{equation}
\rho = \limsup_n \frac{n \ln n}{-\ln |a_n(z)|}
\label{PREL11}
\end{equation}
and (in the case where $0 < \rho < \infty$)
\begin{equation}
\tau = \frac{1}{e \rho}\limsup_n n |a_n(z)|^{\rho/n}
= \frac{e^{\rho - 1}}{\rho}\limsup_n n^{1-\rho} \left|g^{(n)}(z)\right|^{\rho/n},
\label{PREL12}
\end{equation}
respectively, independently of the complex number $z$. An interesting question here is to inquire into the dependence on $z$ of the subsequence(s)
$n_k$ for which the $\limsup$ is atained in \eqref{PREL11} and \eqref{PREL12}.

It is, sometimes, more convenient to write \eqref{PREL11} in the equivalent form
(since $\lim_n |a_n(z)| = 0$ and, hence, $-\ln |a_n(z)|$ is eventually positive)
\begin{equation}
e^{-1/\rho} = \limsup_n |a_n(z)|^{\frac{1}{n \ln n}}
\label{PREL13}
\end{equation}
or, in view of \eqref{PREL10} and the fact that $\lim_n (n!)^{\frac{1}{n \ln n}} = e$,
\begin{equation}
\theta := e^{1 - (1/\rho)} = \limsup_n \left|g^{(n)}(z)\right|^{\frac{1}{n \ln n}}.
\label{PREL14}
\end{equation}
Notice that $\theta = \theta(\rho)$ is smooth and strictly increasing for $\rho \in [0, +\infty]$, with $\theta(0) := \theta(0^+) = 0$ and
$\theta(+\infty) = e$.

Also, if
\begin{equation}
\theta^{\sharp} := \limsup_n \left|\left(g^{\sharp}\right)^{(n)}(z)\right|^{\frac{1}{n \ln n}},
\qquad
z \in \mathbb{C},
\label{PREL14a}
\end{equation}
then we must, obviously, have $\theta^{\sharp} = \theta$ since, as we have seen, $\rho\left(g^{\sharp}\right) = \rho(g) = \rho$. Thus, if we set
\begin{equation}
m_n(r) := \max_{|z| \leq r} \left|g^{(n)}(z)\right|
\qquad
r > 0,
\label{PREL14b}
\end{equation}
then, clearly, $m_n(r) \leq (g^{(n)})^{\sharp}(r) = (g^{\sharp})^{(n)}(r)$ and, therefore
\begin{equation}
\theta = \limsup_n \left[m_n(r)\right]^{\frac{1}{n \ln n}} = \limsup_n \left[\left(g^{\sharp}\right)^{(n)}(r)\right]^{\frac{1}{n \ln n}},
\qquad
r > 0.
\label{PREL14c}
\end{equation}
Likewise, in view of \eqref{PREL12}, since $\tau(g^{\sharp}) = \tau(g) = \tau$ (in the case where $0 < \rho < \infty$) we have
\begin{equation}
\tau = \frac{e^{\rho - 1}}{\rho}\limsup_n n^{1-\rho} \left[m_n(r)\right]^{\rho/n}
= \frac{e^{\rho - 1}}{\rho}\limsup_n n^{1-\rho} \left[\left(g^{\sharp}\right)^{(n)}(r)\right]^{\rho/n},
\quad
r > 0.
\label{PREL14d}
\end{equation}

\subsection{Additional properties of the order}
For the purposes of the present work we need to consider the cases where the $\limsup$ appearing in \eqref{PREL11} is taken over the subsequences
$n = 2k$ and $n = 2k+1$, namely
\begin{equation}
\rho_0(z) := \limsup_k \frac{2k \ln(2k)}{-\ln |a_{2k}(z)|} = 2 \limsup_k \frac{k \ln k}{-\ln |a_{2k}(z)|},
\qquad
z \in \mathbb{C},
\label{PREL15}
\end{equation}
and
\begin{equation}
\rho_1(z) := \limsup_k \frac{(2k+1) \ln(2k+1)}{-\ln |a_{2k+1}(z)|} = 2 \limsup_k \frac{k \ln k}{-\ln |a_{2k+1}(z)|},
\qquad
z \in \mathbb{C},
\label{PREL16}
\end{equation}
respectively. Clearly, for every $z \in \mathbb{C}$ we have
\begin{equation}
\rho = \max\{\rho_0(z), \, \rho_1(z)\}.
\label{PREL17}
\end{equation}

In the spirit of \eqref{PREL14} it is, again, convenient to introduce the quantities
\begin{equation}
\theta_0(z) := \exp\left(1 - \frac{1}{\rho_0(z)}\right) = \limsup_k \left|g^{(2k)}(z)\right|^{\frac{1}{2k \ln k}},
\qquad
z \in \mathbb{C},
\label{PREL18}
\end{equation}
and
\begin{equation}
\theta_1(z) := \exp\left(1 - \frac{1}{\rho_1(z)}\right) = \limsup_k \left|g^{(2k+1)}(z)\right|^{\frac{1}{2k \ln k}},
\qquad
z \in \mathbb{C}.
\label{PREL19}
\end{equation}
Clearly, in view of \eqref{PREL14} we have
\begin{equation}
0 \leq \theta = \max\{\theta_0(z), \, \theta_1(z)\} \leq e.
\label{PREL20}
\end{equation}

We wish to determine how close are the quantities $\rho_0(z)$ and $\rho_1(z)$ to the order $\rho$ of $g(z)$ or, equivalently, how close are
the quantities $\theta_0(z)$ and $\theta_1(z)$ to the constant $\theta$ of \eqref{PREL14}.

Recall that a function $\phi(z)$, defined in a domain $\Omega$ of the complex plane and taking values in $\mathbb{R} \cup \{-\infty\}$, is called
subharmonic (in $\Omega$) if it is locally integrable and for any disk
$D_r(z_0) := \{z \in \mathbb{C} \, : \, |z - z_0| < r\} \subset \Omega$ we have
\begin{equation}
\phi(z_0) \leq \frac{1}{\pi r^2}\int_{D_r(z_0)} \phi(z)\, dx dy
\label{PREL21}
\end{equation}
(here, of course, $z = x + iy = (x, y)$ and the function $\phi(z)$ is subharmonic with respect to the real variables $x$ and $y$). Some authors
require \eqref{PREL21} to hold for almost every $z_0 \in \Omega$, in order to completely characterize subharmonic functions as functions whose
distributinal Laplacian is nonnegative \cite{L-L}. However, such variants of the definition of subharmonicity are nonessential for our analysis.

If $A(z)$ is analytic in a domain $\Omega \subset \mathbb{C}$, then $\ln|A(z)|$ is subharmonic
in $\Omega$ (this follows, e.g., from the facts that (i) $\ln|z - z_0|$ is subharmonic and (ii) if $A(z)$ does not vanish in $\Omega$, then
$\ln|A(z)|$ is harmonic). Also, since $\phi(x) = e^{\alpha x}$, $x \in \mathbb{R}$, is convex for any $\alpha > 0$, Jensen's inequality implies that
$|A(z)|^{\alpha} = e^{\alpha\ln|A(z)|}$ too is subharmonic for any $\alpha > 0$.

\medskip

\textbf{Lemma 1.} The functions $\theta_0(z)$ and $\theta_1(z)$, defined by \eqref{PREL18} and \eqref{PREL19} respectively, are subharmonic in
$\mathbb{C}$.

\smallskip

\textit{Proof}. Let us set
\begin{equation}
\Phi_n(z) := \sup_{k \geq n} \left|g^{(2k)}(z)\right|^{\frac{1}{2k \ln k}},
\qquad
n \geq 2.
\label{PREL22}
\end{equation}
Fix an $r > 0$ and restrict $z \in D_r := D_r(0) = \{z \, : \, |z| \leq r\}$. Then, $|g^{(2k)}(z)| \leq (g^{\sharp})^{(2k)}(r)$. Furthermore, as we have seen,
\begin{equation}
\theta = \limsup_n \left[(g^{\sharp})^{(n)}(r)\right]^{\frac{1}{n \ln n}} \leq e.
\label{PREL23}
\end{equation}
It follows that there is an $M = M(r) > 0$, such that $\Phi_n(z)$ of \eqref{PREL22} is $\leq M$ for all $n \geq 2$ and all $z \in D_r$.

Now, from the discussion preceding Lemma 1 we know that $|g^{(2k)}(z)|^{\frac{1}{2k \ln k}}$ is subharmonic for any $k \geq 2$.
It, then, follows easily that $\Phi_n(z)$ is subharmonic in $D_r$ for every $n \geq 2$ (being finite and the supremum of a sequence of subharmonic
functions). Furthermore, it is obvious that $\Phi_n(z)$ decreases with $n$ and, in view of \eqref{PREL18},
\begin{equation}
\theta_0(z) = \lim_n \Phi_n(z),
\qquad
z \in D_r.
\label{PREL24}
\end{equation}
Therefore, by a simple application of the bounded convergence theorem we can conclude that $\theta_0(z)$ is subharmonic in $D_r$ and, since $r$ is
arbitrary, that $\theta_0(z)$ is subharmonic in $\mathbb{C}$.

In exactly the same way we can show that $\theta_1(z)$ too is subharmonic in $\mathbb{C}$.
\hfill $\blacksquare$

\medskip

\textbf{Remark 1.} (i) Let $\nu := \{n_k\}_{k = 1}^{\infty}$ be a strictly increasing sequence of positive integers and
\begin{equation}
\mathbb{N}_{\nu} := \{n_k \, : \, k \in \mathbb{N}\} = \{n_1, n_2, \ldots\}.
\label{PREL24aa}
\end{equation}
If
\begin{equation}
\theta(z; \nu) := \limsup_k \left|g^{(n_k)}(z)\right|^{\frac{1}{n_k \ln(n_k)}}
= \limsup_{n \in \mathbb{N}_{\nu}} \left|g^{(n)}(z)\right|^{\frac{1}{n\ln(n)}},
\qquad
z \in \mathbb{C},
\label{PREL24a}
\end{equation}
then, by imitating the proof of Lemma 1, we can show that $\theta(z; \nu)$ is subharmonic in $\mathbb{C}$.

(ii) It is a well-known fact \cite{L-L} that a subharmonic function $\phi(z)$ in a domain $\Omega$ is equal to an upper semicontinuous function
$\tilde{\phi}(z)$ for almost every (a.e.) $z \in \Omega$. Therefore, Lemma 1 implies
\begin{equation}
\theta_0(z) = \tilde{\theta}_0(z)
\qquad \text{and} \qquad
\theta_1(z) = \tilde{\theta}_1(z),
\qquad \text{for \;a.e. } z \in \mathbb{C},
\label{PREL25}
\end{equation}
where $\tilde{\theta}_0(z)$ and $\tilde{\theta}_1(z)$ are upper semicontinuous in $\mathbb{C}$.

\medskip

\textbf{Example 1.} Suppose $g(z) = \sin(\lambda z)$, where $\lambda \in \mathbb{C} \setminus \{0\}$. Then, $\rho = 1$ and $\tau = |\lambda|$. Furthermore, since
\begin{equation*}
g^{(2k)}(z) = (-1)^k \lambda^{2k} \sin(\lambda z),
\end{equation*}
we have, in view of \eqref{PREL18},
\begin{equation*}
\theta_0(z)
= \left\{
  \begin{array}{cc}
    1, & \ \lambda z/\pi \in \mathbb{C} \setminus \mathbb{Z}; \\
    0, & \ \lambda z/\pi \in \mathbb{Z}, \\
  \end{array}
\right.
\end{equation*}
where $\mathbb{Z}$ is the set of integers. Obviously, $\theta_0(z)$ is subharmonic and it is equal to $\tilde{\theta}_0(z) \equiv 1$ for all except
for countably many $z \in \mathbb{C}$.

\medskip

We are now ready for a key result.

\medskip

\textbf{Theorem 1.} For an entire function $g(z)$ let $\theta$, $\theta_0(z)$, and $\theta_1(z)$ be as in \eqref{PREL14}, \eqref{PREL18}, and \eqref{PREL19} respectively. Then,
\begin{equation}
\theta_0(z) = \theta
\qquad \text{and} \qquad
\theta_1(z) = \theta
\qquad \text{for \;a.e. } z \in \mathbb{C}.
\label{PREL26}
\end{equation}

\smallskip

\textit{Proof}. By Lemma 1 we have that $\theta_0(z)$ is subharmonic in $\mathbb{C}$. Hence, in view of \eqref{PREL20} we must have
\begin{equation}
\theta_0(w) \leq \frac{1}{\pi r^2} \int_{D_r(w)} \theta_0(z) \, dx dy \leq \theta
\label{PREL27}
\end{equation}
for any $w \in \mathbb{C}$ and any $r > 0$.
Therefore, if for some $w \in \mathbb{C}$ we have that $\theta_0(w) = \theta$, then formula \eqref{PREL27} implies that $\theta_0(z) = \theta$
for a.e. $z \in D_r(w)$, which in turn implies $\theta_0(z) = \theta$ for a.e. $z \in \mathbb{C}$, since $r$ is arbitrary.

More generally, let us only assume that the supremum of $\theta_0(z)$ on some compact subset of $\mathbb{C}$ is $\theta$, namely that there is a
sequence $\{z_n\}_{n=1}^{\infty}$ with $\lim_n z_n = z_{\ast} \in \mathbb{C}$ and $\lim_n \theta_0(z_n) = \theta$. We will show that, we must, again,
have $\theta_0(z) = \theta$ for a.e. $z \in \mathbb{C}$.

Fix a disk $D_r(z_{\ast})$ and consider the disks $D_n := D_{r_n}(z_n)$, $n = 1, 2, \ldots$, so that $r_n$ is the largest radius satisfying
$D_n \subset D_r(z_{\ast})$. Using $w = z_n$ and $D_r(w) = D_n$ in \eqref{PREL27} yields
\begin{equation}
\theta_0(z_n) \leq \frac{1}{\pi r_n^2} \int_{D_n} \theta_0(z) \, dx dy
\leq \frac{1}{\pi r_n^2} \int_{D_r(z_{\ast})}  \theta_0(z) \, dx dy
\leq \frac{r^2}{r_n^2} \, \theta,
\qquad
n \geq 1,
\label{PREL28}
\end{equation}
thus, by letting $n \to \infty$ we obtain
\begin{equation}
\theta \leq \frac{1}{\pi r^2} \int_{D_r(z_{\ast})} \theta_0(z) \, dx dy
\leq \theta,
\label{PREL29}
\end{equation}
which tells us that $\theta_0(z) = \theta$ for a.e. $z \in D_r(z_{\ast})$ and, consequently, that $\theta_0(z) = \theta$ for a.e. $z \in \mathbb{C}$.

In the same manner we can show that if the supremum of $\theta_1(z)$ on some compact subset of $\mathbb{C}$ is $\theta$, then
$\theta_1(z) = \theta$ for a.e. $z \in \mathbb{C}$.

Finally, we will show that the assumption
\begin{equation}
\Theta_0(r) := \sup_{|z| \leq r} \theta_0(z) < \theta
\qquad \text{for every }\; r > 0
\label{PREL30}
\end{equation}
leads to a contradiction.

For a given $r > 0$ let us assume \eqref{PREL30} and fix an $\varepsilon > 0$ so that
\begin{equation}
\Theta_0(r) + \varepsilon < \theta.
\label{PREL31}
\end{equation}
Then, by \eqref{PREL18} and \eqref{PREL30} we get that for every $z \in D_r(0)$ there is an integer $K = K(z)$ such that
\begin{equation}
\sup_{k \geq K(z)} \left|g^{(2k)}(z)\right|^{\frac{1}{2k \ln k}} < \Theta_0(r) + \varepsilon.
\label{PREL32}
\end{equation}
It follows that if we set
\begin{equation}
G_m := \left\{z \in D_r(0) \, : \, \sup_{k \geq m} \left|g^{(2k)}(z)\right|^{\frac{1}{2k \ln k}} < \Theta_0(r) + \varepsilon\right\},
\label{PREL33}
\end{equation}
then
\begin{equation}
\bigcup_{m=1}^{\infty} G_m = D_r(0).
\label{PREL34}
\end{equation}
Thus, there is a $m_0$ for which the set $G_{m_0}$ has positive (Lebesgue) measure.

Now, as we have seen in the proof of Lemma 1, the function
\begin{equation}
\phi(z) :=  \sup_{k \geq m_0} \left|g^{(2k)}(z)\right|^{\frac{1}{2k \ln k}}
\label{PREL35}
\end{equation}
is subharmonic. Thus, as we have mentioned in Remark 1(ii) there is an upper semicontinuous function $\tilde{\phi}(z)$ such that
$\phi(z) = \tilde{\phi}(z)$ for a.e. $z \in D_r(0)$. Therefore, the sets
\begin{equation}
G_{m_0} = \left\{z \in D_r(0) \, : \, \phi(z) < \Theta_0(r) + \varepsilon\right\}
\quad \text{and} \quad
\tilde{G} := \left\{z \in D_r(0) \, : \, \tilde{\phi}(z) < \Theta_0(r) + \varepsilon\right\}
\label{PREL36}
\end{equation}
differ by a set of measure $0$, i.e. the set $G_{m_0} \vartriangle \tilde{G}$ has zero (Lebesgue) measure. Furthermore, the upper semicontinuity
of $\tilde{\phi}(z)$ implies \cite{L-L} that $\tilde{G}$ is open (and nonempty since $G_{m_0}$ has positive measure). Therefore, any open disk
$D_{\delta}(z_0) \subset \tilde{G}$ lies almost entirely in $G_{m_0}$ in the sense that their symmetric difference has measure $0$ (in other words,
the area of $D_{\delta}(z_0) \cap G_{m_0}$ is equal to the area of $D_{\delta}(z_0)$, namely $\pi\delta^2$).

We continue by noticing that the assumption \eqref{PREL30} implies that $\theta_0(z) < \theta$ for all $z \in \mathbb{C}$ and, hence, we must have
\begin{equation}
\theta_1(z) \equiv \theta.
\label{PREL31a}
\end{equation}
Let $D_{\delta}(z_0)$ be a disk  (with $\delta > 0$) such that $D_{\delta}(z_0) \subset \tilde{G}$. If $\Gamma$ is the boundary of $D_{\delta}(z_0)$,
then due to the previous discussion we can arrange it so that the symmetric difference of $G_{m_0}$ and $\Gamma$ has one-dimensional measure $0$ (in
other words, the ``$\,$length" (i.e. the one-dimensional measure) of $\Gamma \cap G_{m_0}$ is equal to the length of $\Gamma$, namely $2\pi\delta$).

Now, by Cauchy's integral formula we have
\begin{equation}
g^{(2k+1)}(z_0) = \frac{1}{2\pi i} \oint_{\Gamma} \frac{g^{(2k)}(z)}{(z - z_0)^2} \, dz,
\qquad
k \geq 0,
\label{PREL32a}
\end{equation}
Taking absolute values in \eqref{PREL32a} yields
\begin{equation}
\left|g^{(2k+1)}(z_0)\right| \leq \frac{1}{2\pi \delta^2} \oint_{\Gamma} \left|g^{(2k)}(z)\right| \, ds,
\qquad
k \geq 0,
\label{PREL33a}
\end{equation}
where $ds$ is the arc-length element of $\Gamma$.

By \eqref{PREL33} we have
\begin{equation}
\left|g^{(2k)}(z)\right| < \left[\Theta_0(r) + \varepsilon\right]^{2k \ln k}
\label{PREL37}
\end{equation}
for all $k \geq m_0$ and a.e. $z \in \Gamma$. Thus, by using \eqref{PREL37} in in \eqref{PREL33a} we obtain
\begin{equation}
\left|g^{(2k+1)}(z_0)\right| \leq \frac{1}{2\pi \delta^2} \oint_{\Gamma} \left[\Theta_0(r) + \varepsilon\right]^{2k \ln k} \, ds
= \frac{\left[\Theta_0(r) + \varepsilon\right]^{2k \ln k}}{\delta},
\quad
k \geq m_0,
\label{PREL38}
\end{equation}
or
\begin{equation}
\left|g^{(2k+1)}(z_0)\right|^{\frac{1}{2k \ln k}} \leq \left(\frac{1}{\delta}\right)^{\frac{1}{2k \ln k}}\left[\Theta_0(r) + \varepsilon\right],
\qquad
k \geq m_0
\label{PREL39}
\end{equation}
Therefore,
\begin{equation}
\theta_1(z_0) = \limsup_k \left|g^{(2k+1)}(z_0)\right|^{\frac{1}{2k \ln k}} \leq \Theta_0(r) + \varepsilon < \theta,
\label{PREL40}
\end{equation}
which contradicts \eqref{PREL31a}. Hence, the assumption \eqref{PREL30} is false and we must have
$\sup_{|z| \leq r} \theta_0(z) = \theta$ for some $r > 0$, which, as we have seen earlier in the proof, implies $\theta_0(z) = \theta$
for a.e. $z \in \mathbb{C}$.

The proof that $\theta_1(z) = \theta$ for a.e. $z \in \mathbb{C}$ is exactly the same.
\hfill $\blacksquare$

\medskip

\textbf{Remark 2.} In view of \eqref{PREL18} and \eqref{PREL19} Theorem 1 implies immediately that if
\begin{equation}
\mathcal{F}_0 := \{z \in \mathbb{C} \,:\, \rho_0(z) = \rho\}
\qquad \text{and} \qquad
\mathcal{F}_1 := \{z \in \mathbb{C} \,:\, \rho_1(z) = \rho\},
\label{PREL41}
\end{equation}
where $\rho$ is the order of $g(z)$ (thus $0 \leq \rho \leq \infty$) and the quantities $\rho_0(z)$ and $\rho_1(z)$ are defined in \eqref{PREL15} and \eqref{PREL16} respectively, then both sets $\mathcal{F}_0$ and $\mathcal{F}_1$ have full measure (and by formula \eqref{PREL17} we have
$\mathcal{F}_0 \cup \mathcal{F}_1 = \mathbb{C}$); in other words the sets $\mathcal{F}_0^c := \mathbb{C} \setminus \mathcal{F}_0$
and $\mathcal{F}_1^c := \mathbb{C} \setminus \mathcal{F}_1$ have Lebesque measure (i.e. area) zero.

\medskip

\textbf{Open Question 1.} Are the sets $\mathcal{F}_0^c$ and $\mathcal{F}_1^c$ nowhere dense in $\mathbb{C}$? Are they countable?

\medskip

\textbf{Remark 3.} Let $\theta(z; \nu)$ be as in \eqref{PREL24a} (see Remark 1). Then, under some general assumptions on the sequence
$\nu := \{n_k\}_{k = 1}^{\infty}$ the proof of Theorem 1 remains valid and we will have $\theta(z; \nu) = \theta$ for a.e. $z \in \mathbb{C}$.
For example, if
\begin{equation*}
\frac{n_{k-1}}{n_k} \rightarrow 1
\qquad \text{as }\;
k \to \infty,
\end{equation*}
then $\theta(z; \nu) = \theta$ for a.e. $z \in \mathbb{C}$.

\subsection{Additional properties of the type}
We, now, turn our attention to the type of $g(z)$. Of course, we need to assume that $0 < \rho < \infty$.

In view of \eqref{PREL12} and inspired by \eqref{PREL15} and \eqref{PREL16} we set
\begin{align}
\tau_0(z) :&= \frac{1}{e \rho}\limsup_k 2k \, |a_{2k}(z)|^{\rho/2k}
= \frac{2}{e \rho}\limsup_k k \, |a_{2k}(z)|^{\rho/2k}
\nonumber
\\
&= \frac{(e/2)^{\rho - 1}}{\rho}\limsup_k k^{1-\rho} \left|g^{(2k)}(z)\right|^{\rho/2k},
\qquad
z \in \mathbb{C}
\label{PREL42}
\end{align}
and
\begin{align}
\tau_1(z) :&= \frac{1}{e \rho}\limsup_k (2k+1) \, |a_{2k+1}(z)|^{\frac{\rho}{2k+1}}
= \frac{2}{e \rho}\limsup_k k \, |a_{2k+1}(z)|^{\frac{\rho}{2k}}
\nonumber
\\
&= \frac{(e/2)^{\rho - 1}}{\rho}\limsup_k k^{1-\rho} \left|g^{(2k+1)}(z)\right|^{\rho/2k},
\qquad
z \in \mathbb{C},
\label{PREL43}
\end{align}
where $a_n(z)$ is given by \eqref{PREL10}. From \eqref{PREL12} it is obvious that
\begin{equation}
\tau = \max\{\tau_0(z), \, \tau_1(z)\}.
\label{PREL44}
\end{equation}

The following theorem gives a property of the type of $g(z)$ which is the analog of the property regarding the order of $g(z)$ established in
Theorem 1 (see also Remark 2).

\medskip

\textbf{Theorem 2.} Let $\tau$ be the type of the entire function $g(z)$, while $\tau_0(z)$ and $\tau_1(z)$ be as in \eqref{PREL42} and \eqref{PREL43} respectively.

(i) If $\tau < \infty$, then
\begin{equation}
\tau_0(z) = \tau
\qquad \text{and} \qquad
\tau_1(z) = \tau
\qquad \text{for \;a.e. } z \in \mathbb{C}.
\label{PREL45}
\end{equation}

(ii) If $\tau = \infty$, then there exists a dense $G_{\delta}$ (therefore uncountable) subset $\mathcal{U}_{\infty}$ of $\mathbb{C}$ such that
\begin{equation}
\tau_0(z) = \infty
\qquad \text{and} \qquad
\tau_1(z) = \infty
\qquad \text{for all } z \in \mathcal{U}_{\infty}.
\label{PREL45a}
\end{equation}

\smallskip

\textit{Proof}. (ii) If $\tau < \infty$, we can follow the proof of Lemma 1 in order to show that $\tau_0(z)$ and $\tau_1(z)$ are subharmonic in $\mathbb{C}$. Then, by imitating the proof of Theorem 1 we can easily obtain \eqref{PREL45}.

(ii) Suppose $\tau = \infty$. Then, in view of \eqref{PREL12} we have
\begin{equation}
\sigma(z) := \sup_n n^{1-\rho} \left|g^{(n)}(z)\right|^{\rho/n} \equiv \infty.
\label{2PROOF1}
\end{equation}
We introduce the quantities
\begin{equation}
\sigma_0(z) := \sup_k \, (2k)^{1-\rho} \left|g^{(2k)}(z)\right|^{\rho/2k},
\qquad
z \in \mathbb{C}
\label{2PROOF2}
\end{equation}
and
\begin{equation}
\sigma_1(z) := \sup_k \, (2k+1)^{1-\rho} \left|g^{(2k+1)}(z)\right|^{\rho/(2k+1)},
\qquad
z \in \mathbb{C},
\label{2PROOF3}
\end{equation}
so that
\begin{equation}
\max\{\sigma_0(z), \, \sigma_1(z)\} \equiv \infty.
\label{2PROOF4}
\end{equation}
Observe that for $j = 0,1$ we have $\tau_j(z) = \infty$ if and only if $\sigma_j(z) = \infty$. Therefore, it suffices to prove \eqref{PREL45a} for
$\sigma_j(z)$ in place of $\tau_j(z)$.

Suppose that for some disk $D$ we had
\begin{equation}
\sup_{z \in D} \sigma_0(z) < \infty.
\label{2PROOF6}
\end{equation}
Then, \eqref{2PROOF4} would imply that $\sigma_1(z) = \infty$ for all $z \in D$, which, it can be shown to be impossible under \eqref{2PROOF6} by
following the approach used in the proof of Theorem 1, starting with formula \eqref{PREL32a}. Therefore,
\begin{equation}
\sup_{z \in D} \sigma_0(z) = \infty
\qquad \text{for any disk }\; D.
\label{2PROOF7}
\end{equation}
Now, formula \eqref{2PROOF2} implies that $\sigma_0(z)$ is lower semicontinuous on $\mathbb{C}$ (being the supremum of continuous functions). Hence,
the set
\begin{equation}
G_N := \{z \in \mathbb{C} \,:\, \sigma_0(z) > N\}.
\label{2PROOF9}
\end{equation}
is open. Furthermore, by \eqref{2PROOF7} we have that $G_N$ is dense in $\mathbb{C}$ and, therefore, the set
\begin{equation}
\{z \in \mathbb{C} \,:\, \sigma_0(z) = \infty\} = \bigcap_{N=1}^{\infty} G_N
\label{2PROOF10}
\end{equation}
is a dense $G_{\delta}$ subset of $\mathbb{C}$. Likewise, $\{z \in \mathbb{C} \,:\, \sigma_1(z) = \infty\}$ is dense and $G_{\delta}$.
Since the intersection of two dense $G_{\delta}$ sets is again a dense $G_{\delta}$ set, we can satisfy \eqref{PREL45a} by taking
\begin{equation}
\mathcal{U}_{\infty} = \{z \in \mathbb{C} \,:\, \sigma_0(z) = \infty\} \cap \{z \in \mathbb{C} \,:\, \sigma_1(z) = \infty\}.
\label{2PROOF11}
\end{equation}
\hfill $\blacksquare$

\medskip

\textbf{Remark 4.} As we have seen, the functions $\left|g^{(n)}(z)\right|^{\rho/n}$, $n = 1, 2, \ldots$, are subharmonic. It follows that $\sigma_j(z)$, $j = 0,1$,
satisfy \eqref{PREL21}, namely
\begin{equation}
\sigma_j(z_0) \leq \frac{1}{\pi r^2}\int_{D_r(z_0)} \sigma_j(z)\, dx dy
\label{2PROOF5}
\end{equation}
for any disk $D_r(z_0)$. Notice, however, that $\sigma_0(z)$ or $\sigma_1(z)$ may not be subharmonic, since they may become infinite for some $z$ or
they may not be locally integrable.
By using \eqref{2PROOF7} in \eqref{2PROOF5}, and aruing as in the beginning of the proof of Theorem 1, we can conclude that
\begin{equation}
\int_D \sigma_j(z)\, dx dy = \infty
\qquad \text{for any disk }\; D
\label{2PROOF8}
\end{equation}
for $j = 0, 1$.
Actually, with the help of Poisson integral formula for harmonic functions (and the fact that in any sufficiently smooth domain a sunharmonic function is dominated by the harmonic function with the same boundary values) we can get a stronger version of \eqref{2PROOF8}, namely
\begin{equation}
\int_{\Gamma} \sigma_j(z)\, ds = \infty
\qquad \text{for any circle }\; \Gamma.
\label{2PROOF8a}
\end{equation}
Furthermore, since $\ln |g^{(n)}(z)|$ is subharmonic, we can work with $\ln \sigma_j(z)$ instead of $\sigma_j(z)$ and conclude that
\begin{equation}
\int_{\Gamma} \ln\sigma_j(z)\, ds = \infty
\qquad \text{for any circle }\; \Gamma
\label{2PROOF8b}
\end{equation}
for $j = 0, 1$. However, in spite of \eqref{2PROOF8b}, the question whether $\tau_0(z) = \infty$ and $\tau_1(z) = \infty$ for a.e. $z \in \mathbb{C}$
remains open in the case where $\tau = \infty$.

\medskip

\subsection{Canonical products}
In this short subsection we review some basic facts regarding canonical products.

Suppose $z_1, z_2, \ldots$ is a finite or infinite sequence of non-zero complex numbers such that
\begin{equation}
\sigma := \inf\left\{s \geq 0 \, : \, \sum_{k \geq 1} \frac{1}{|z_k|^s} < \infty\right\} < \infty.
\label{RO1}
\end{equation}
Then, the \emph{canonical product} (\emph{associated to} $\{z_k\}_{k \geq 1}$) is the quantity
\begin{equation}
\Pi(z) := \prod_{k \geq 1} e_p\left(\frac{z}{z_k}\right),
\label{RO2}
\end{equation}
where
\begin{equation}
e_0(z) := 1 - z,
\qquad
e_p(z) := (1 - z) \exp\left(\frac{z}{1} + \frac{z^2}{2} + \cdots + \frac{z^p}{p}\right),
\quad
p > 0,
\label{RO3}
\end{equation}
and $p$ is related to $\sigma$ as follows:

(i) If $\sigma$ is not an integer, then $p = \lfloor\sigma\rfloor$.

(ii) If $\sigma$ is an integer and
\begin{equation}
\sum_{k \geq 1} \frac{1}{|z_k|^{\sigma}} = \infty,
\label{RO4a}
\end{equation}
then $p = \sigma$.

(iii) If $\sigma$ is an integer and
\begin{equation}
\sum_{k \geq 1} \frac{1}{|z_k|^{\sigma}} < \infty,
\label{RO4b}
\end{equation}
then $p = \max\{\sigma - 1, 0\}$ (notice that, under \eqref{RO4b}, $\sigma = 0$ if and only if $\{z_k\}$ is a finite sequence).

It is a well-known fact in complex analysis \cite{H} that the canonical product of \eqref{RO2} is entire in $z$ of order $\sigma$ and, furthermore, in the case $\sigma > 0$, this entire function is of minimal type if $\sum_{k \geq 1} |z_k|^{-\sigma} < \infty$ (the converse may not be true).

\medskip

\textbf{Remark 5.} In view of \eqref{RO3} formula \eqref{RO2} can be written as
\begin{equation}
\Pi(z) = \prod_{k \geq 1} \left(1 - \frac{z}{z_k}\right) \exp\left(\frac{z}{z_k} + \frac{z^2}{2 z_k^2} + \cdots + \frac{z^p}{p z_k^p}\right).
\label{RO2a}
\end{equation}
Obviously,
\begin{equation}
\Pi(0) = 1.
\label{RO2b}
\end{equation}
Also, in case where $p \geq 1$,
\begin{equation}
\frac{\Pi'(z)}{\Pi(z)} = \sum_{k \geq 1} \left(\frac{1}{z - z_k} + \frac{1}{z_k} + \frac{z}{z_k^2} + \cdots + \frac{z^{p-1}}{z_k^p}\right).
\label{RO2c}
\end{equation}
Thus
\begin{equation}
\frac{\Pi'(0)}{\Pi(0)} = 0,
\qquad\text{which implies}\qquad
\Pi'(0) = 0
\label{RO2d}
\end{equation}
and, more generally,
\begin{equation*}
\left.\frac{d^{r-1}}{dz^{r-1}}\left[\frac{\Pi'(z)}{\Pi(z)}\right]\right|_{z=0} = 0
\qquad \text{for }\; r = 1, 2, \ldots, p,
\end{equation*}
which implies
\begin{equation}
\Pi^{(r)}(z) = 0
\qquad \text{for }\; r = 1, 2, \ldots, p.
\label{RO2f}
\end{equation}
Finally, let us notice that if
\begin{equation}
g(z) := e^{A_1 z + A_2 z^2 + \cdots + A_m z^m} \Pi(z),
\label{RO2g}
\end{equation}
where $A_1, A_2, \ldots, A_m$ are complex constants and $\Pi(z)$ is as in \eqref{RO2a} with $p \geq 1$, then
\begin{equation}
g(0) = 1
\qquad \text{and} \qquad g'(0) = A_1
\label{RO2gg}
\end{equation}
(if $p \geq 2$, then $g''(0) = 2A_2 + A_1^2$).

\medskip

\subsection{Some notation and terminology}
From now on we will use the following notation$/$terminology for typographical convenience.

\medskip

Let $g(z)$ be an entire function of order $\rho \in (0, \infty)$ and type $\tau$.

(i) If $\tau = 0$ we will say that the \emph{exact order} of $g(z)$ is $\rho^-$.

(ii) If $0 < \tau < \infty$ we will say that the \emph{exact order} of $g(z)$ is $\rho$.

(iii) If $\tau = \infty$ we will say that the \emph{exact order} of $g(z)$ is $\rho^+$.

\medskip

In the case where $\rho = 0$ or $\rho = \infty$ we can consider of the notion of exact order as equivalent to the order.

For instance, the statement that the exact order of $g(z)$ is in $[0, 2^-]$ means that the order is between $0$ and $2$ (included), and in the case
where it is equal to $2$ the type of $g(z)$ is $0$.

\section{Order and type considerations for entire caloric functions}

\subsection{The $t$- and the $z$-power series of a caloric function}
For $t_0 = 0$ the second equality in \eqref{P2} can be written in the (Taylor) form
\begin{equation}
F(t, z) = \sum_{j \geq 0} \frac{\partial_t^j F(0, z)}{j!} \, t^j.
\label{O1}
\end{equation}
If $F(t, z)$ satisfies the heat equation \eqref{P3}, we can use \eqref{P4} in \eqref{O1} and get the following expansion in powers of $t$
\begin{equation}
F(t, z) = \sum_{j \geq 0} \frac{\partial_z^{2j} F(0, z)}{j!} \, t^j = \sum_{j \geq 0} \frac{f^{(2j)}(z)}{j!} \, t^j,
\label{O2}
\end{equation}
where (recall \eqref{P3a}),
\begin{equation}
f(z) = F(0, z) = \sum_{k \geq 0} c_k z^k.
\label{O3}
\end{equation}

Let us also expand $F(t, z)$ in powers of $z$. For $z_0 = 0$ the first equality in \eqref{P2} can be written in the
(Taylor) form
\begin{equation}
F(t, z) = \sum_{k \geq 0} \frac{\partial_z^k F(t, 0)}{k!} \, z^k
= \sum_{k \geq 0} \frac{\partial_z^{2k} F(t, 0)}{(2k)!} \, z^{2k}
+ \sum_{k \geq 0} \frac{\partial_z^{2k} \left[\partial_z F(t, 0)\right]}{(2k+1)!} \, z^{2k+1}.
\label{TO1}
\end{equation}
If $F(t, z)$ is caloric, so is $\partial_z F(t, z)$. Thus, we can use \eqref{P4} in \eqref{TO1} and get the following expansion
\begin{equation}
F(t, z)
= \sum_{k \geq 0} \frac{\partial_t^k F(t, 0)}{(2k)!} \, z^{2k}
+ \sum_{k \geq 0} \frac{\partial_t^k \left[\partial_z F(t, 0)\right]}{(2k+1)!} \, z^{2k+1},
\label{TO2}
\end{equation}
or
\begin{equation}
F(t, z)
= \sum_{k \geq 0} \frac{\phi^{(k)}(t)}{(2k)!} \, z^{2k}
+ \sum_{k \geq 0} \frac{\psi^{(k)}(t)}{(2k+1)!} \, z^{2k+1},
\label{TO3}
\end{equation}
where we have set
\begin{equation}
\phi(t) := F(t, 0)
\qquad \text{and} \qquad
\psi(t) := \partial_z F(t, 0).
\label{TO4}
\end{equation}
We can say that formula \eqref{TO3} is the counterpart of \eqref{O2}.

\medskip

\textbf{Remark 6.} If $F(t, z)$ satisfies the heat equation \eqref{P3}, so does $F(t, -z)$. Consequently, the $z$-even and $z$-odd parts of $F(t, z)$, namely
\begin{equation}
F_e(t, z) := \frac{F(t, z) + F(t, -z)}{2}
\qquad \text{and} \qquad
F_o(t, z) := \frac{F(t, z) - F(t, -z)}{2}
\label{TO5}
\end{equation}
respectively, also satisfy \eqref{P3}, and, moreover, from \eqref{TO3} we get immediately that
\begin{equation}
F_e(t, z) = \sum_{k \geq 0} \frac{\phi^{(k)}(t)}{(2k)!} \, z^{2k}
\qquad \text{and} \qquad
F_o(t, z) = \sum_{k \geq 0} \frac{\psi^{(k)}(t)}{(2k+1)!} \, z^{2k+1}.
\label{TO6}
\end{equation}

\medskip

Finally, let us notice that if $\phi(t)$ and $\psi(t)$ are two arbitrary entire functions, then, for any fixed $t \in \mathbb{C}$ the radii of convergence of the two $z$-power series in \eqref{TO3} are infinite (and the same is true for the $t$-derivatives of those $z$-power series), hence
both series are entire in $(t, z)$ and satisfy
the heat equation \eqref{P3}. Thus, every entire solution $F(t, z)$ of \eqref{P3} is determined uniquely by a pair of
(arbitrary) entire functions $\phi(t)$ and $\psi(t)$ via \eqref{TO3}-\eqref{TO4}. As we will see in the next subsection, the
situation is quite different if we want to determine an entire caloric solution $F(t, z)$ from the
initial condition $f(z)$ via \eqref{O2}, since in this case the entire function $f(z)$ cannot be arbitrary.

\subsection{The order of the initial condition $f(z)$}
Recall that in view of \eqref{PREL8} the order $\rho$ of the entire function $f(z)$ can be expressed as
\begin{equation}
\rho = \limsup_k \frac{k \, \ln k}{-\ln|c_k|},
\label{O4}
\end{equation}
where, in view of \eqref{O3},
\begin{equation}
c_k = \frac{f^{(k)}(0)}{k!},
\qquad
k \geq 0.
\label{O4a}
\end{equation}

\medskip

\textbf{Theorem 3.} Suppose $f(z)$ is the initial condition of an entire solution $F(t, z)$ of the heat equation \eqref{P3}. Then the order $\rho$
of $f(z)$ satisfies $0 \leq \rho \leq 2$. Furthermore, in the extreme case $\rho = 2$ the type $\tau$ of $f(z)$ must be minimal, i.e. zero. Thus,
under the terminology introduced in Subsection 2.4, the exact order of $f(z)$ is in $[0, 2^-]$.

Conversely, if $f(z)$ is an entire function whose exact order is in $[0, 2^-]$, then $F(t, z)$ given by \eqref{O2} is entire in
$(t, z)$ and satisfies the heat equation \eqref{P3} with $F(0, z) = f(z)$.

\smallskip

\textit{Proof}. We need to consider the cases $k = 2j$ and $k = 2j+1$ separately.

Observe that $\rho = \max\{\rho_0, \rho_1\}$ where
\begin{equation}
\rho_0 := \limsup_j \frac{2j \, \ln (2j)}{-\ln |c_{2j}|} = 2 \, \limsup_j \frac{j \, \ln j}{-\ln |c_{2j}|}
\label{O4b}
\end{equation}
and
\begin{equation}
\rho_1 := \limsup_j \frac{(2j+1) \, \ln (2j+1)}{-\ln |c_{2j+1}|} = 2 \, \limsup_j \frac{j \, \ln j}{-\ln |c_{2j+1}|}.
\label{O4c}
\end{equation}
Let us first estimate $\rho_0$. For $z = 0$ formula \eqref{O2}, in view of \eqref{O4a}, becomes
\begin{equation}
F(t, 0) = \sum_{j \geq 0} \frac{t^j}{j!} \, f^{(2j)}(0) = \sum_{j \geq 0} \frac{(2j)! \, c_{2j}}{j!} \, t^j.
\label{O5}
\end{equation}
Since the second power series in \eqref{O5} converges for all $t \in \mathbb{C}$, i.e. has an infinite radius of convergence, we must have
\begin{equation}
\limsup_j \left|\frac{(2j)! \, c_{2j}}{j!}\right|^{1/j} = 0
\qquad
\left(\text{hence,}\quad \lim_j \left|\frac{(2j)! \, c_{2j}}{j!}\right|^{1/j} = 0\right)
\label{O6}
\end{equation}
equivalently, there is a sequence $\varepsilon_j \to 0$ such that
\begin{equation}
|c_{2j}|^{1/j} < \varepsilon_j \left[\frac{j!}{(2j)!}\right]^{1/j}.
\label{O6a}
\end{equation}
By applying Stirling's asymptotic formula for $n!$ to \eqref{O6a} we get that there is a sequence $\varepsilon_j \to 0$ (not necessarily the same
$\varepsilon_j$ appearing in \eqref{O6a}) such that
\begin{equation}
|c_{2j}|^{1/j} < \frac{\varepsilon_j}{j}.
\label{O6b}
\end{equation}
Now, \eqref{O4b} tells us that $2 / \rho_0$ is the supremum of all exponents $r$ such that
\begin{equation}
|c_{2j}|^{1/j} < \frac{1}{j^r}
\qquad \text{for every sufficiently large }\; j.
\label{O7}
\end{equation}
Therefore, by comparing \eqref{O6b} and \eqref{O7} we get that
\begin{equation}
2 / \rho_0 \geq 1,
\qquad \text{i.e.} \qquad
\rho_0 \leq 2.
\label{O8}
\end{equation}
To determine $\rho_1$ we differentiate \eqref{O2} with respect to $z$ and arrive at
\begin{equation}
\partial_z F(t, z) = \sum_{j \geq 0} \frac{t^j}{j!} \, f^{(2j+1)}(z).
\label{O9}
\end{equation}
Then, in exactly the same way as in the case of $\rho_0$ we conclude that there is a sequence $\varepsilon_j \to 0$ such that
\begin{equation}
|c_{2j+1}|^{1/j} < \frac{\varepsilon_j}{j}.
\label{O9a}
\end{equation}
Hence
\begin{equation}
\rho_1 \leq 2
\label{O10}
\end{equation}
and, therefore, the order of $f(z)$ satisfies $\rho \leq 2$.

Next, let us estimate the type $\tau$ of $f(z)$ in the case where $\rho = 2$. Recall that if $\rho = 2$, then
\begin{equation}
\tau = \frac{1}{2e} \, \limsup_k k \, |c_k|^{2/k}.
\label{O11}
\end{equation}
As in the case of $\rho$ we, again, need to consider the cases $k = 2j$ and $k = 2j+1$ odd separately. We, thus, observe that
$\tau = \max\{\tau_0, \tau_1\}$ where
\begin{equation}
\tau_0 := \frac{1}{e} \, \limsup_j j \, |c_{2j}|^{1/j}
\qquad \text{and} \qquad
\tau_1 := \frac{1}{e} \, \limsup_j j \, |c_{2j+1}|^{1/j}.
\label{O12}
\end{equation}
Then, by \eqref{O6b} and \eqref{O9a} we get immediately that $\tau_0 = \tau_1 = 0$ and, consequently $\tau = 0$. In other words, if $f(z)$ is of
order $2$, then it is of minimal type.

The second part of the theorem is easy (see also Remark 8 below).
\hfill $\blacksquare$

\medskip

Using the Hadamard Factorization Theorem \cite{H} we get immediately the following corollary of Theorem 3.

\medskip

\textbf{Corollary 1.} Suppose $f(z) \not\equiv 0$ is the initial condition of an entire solution $F(t, z)$ of the heat equation \eqref{P3}. Then
$f(z)$ has the
form
\begin{equation}
f(z) = e^{\lambda z + \beta} z^d \, \Pi(z),
\label{CO1}
\end{equation}
where $\lambda$ and $\beta$ are complex constants, $d$ is a nonnegative integer and $\Pi(z)$ is a canonical product whose exact order is in
$[0, 2^-]$.

\medskip

Finally notice that from Subsection 2.1 it follows that the canonical product $\Pi(z)$ of \eqref{CO1} must be of the form
\begin{equation}
\Pi(z) = \prod_{k \geq 1} e_p\left(\frac{z}{a_k}\right),
\qquad \text{where }\;
p = 0 \ \text{or } 1 \ \text{or } 2,
\label{CO1a}
\end{equation}
and the case $p = 2$ can happen only if the exact order of $\Pi(z)$ is $2^-$ and $\sum_{k \geq 1} |z_k|^{-2} = \infty$.

\subsection{The order and type of a caloric function $F(t, z)$ when $t$ or $z$ is fixed}
Suppose we freeze $t \in \mathbb{C}$ and we consider
\begin{equation}
\rho_z := \text{ord}_z F(t, z)
\qquad\text{and}\qquad
\tau_z := \text{type}_z F(t, z),
\label{IO1}
\end{equation}
namely the order and the type of $F(t, \cdot)$, which for convenience we call the \emph{$z$-order} and \emph{$z$-type of}
$F(t, z)$. One might think that these quantities depend on $t$. However, by applying
Theorem 1.1 of \cite{C-K-K} for the choice $f(D) = e^{tD^2}$ (in the notation of Theorem 1.1 of \cite{C-K-K}) and noticing that
$F(t, z) = e^{tD^2} F(0, z)$ with $D = \partial_z$, we get that $\rho_z$ and $\tau_z$ are, actually, independent of $t$ (actually, this also follows
from \eqref{TO3}). We can, therefore, choose $t = 0$ and conclude that
\begin{equation}
\rho_z = \rho
\qquad\text{and}\qquad
\tau_z = \tau,
\label{IO1a}
\end{equation}
where $\rho$ and $\tau$ are the order and type of $f(z) = F(0, z)$ respectively. And then, Theorem 3 implies immediately that the exact $z$-order
of $F(t, z)$ is in $[0, 2^-]$.

The following open question envisages a refinement of formula \eqref{IO1a}.

\medskip

\textbf{Open Question 2.} Let $z_1(t), z_2(t), \ldots$  be the zeros of $F(t, z)$ (viewed as an entire function of $z$). Is it true that for any
$t_1, t_2 \in \mathbb{C}$ and any $\alpha > 0$ we have that
\begin{equation*}
\sum_{k \geq 1}\frac{1}{|z_k(t_1)|^{\alpha}} < \infty
\qquad\text{if and only if }\qquad
\sum_{k \geq 1}\frac{1}{|z_k(t_2)|^{\alpha}} < \infty?
\end{equation*}
More specifically, if the zeros are infinitely many and arranged so that $|z_k(t)| \leq |z_{k+1}(t)$ for all $k \geq 1$,
how close are the asymptotics of the sequences $\{z_k(t_1)\}_{k \geq 1}$ and $\{z_k(t_2)\}_{k \geq 1}$ as $k \to \infty$?

\medskip

Regarding the order and the type of $F(t, z)$ viewed as a function of $t$
things are more complicated, since these quantities may depend on $z$. For instance it is clear from \eqref{TO3}-\eqref{TO4} that there are nontrivial caloric functions $F(t, z)$ such that, for a given $z_0$ the quantity $F(t, z_0)$ can vanish for all $t \in \mathbb{C}$ (e.g., just think of
$F(t, z) = e^{-\lambda^2 t} \sin(\lambda z)$, for $z_0 = k\pi/\lambda$, $k \in \mathbb{Z}$).

The goal of this subsection is to clarify the notions of $t$-order and $t$-type of an entire caloric function $F(t, z)$ and to relate them to
its $z$-order and its $z$-type respectively.

\medskip

\textbf{Theorem 4.} Let $F(t, z)$ be an entire solution of the heat equation \eqref{P3}. For a fixed $z \in \mathbb{C}$ we consider the orders of
$F(t, z)$ and $\partial_z F(t, z)$ viewed as functions of $t$, namely
\begin{equation}
\rho_{t,0}(z) := \text{ord}_t F(t, z),
\qquad
\rho_{t,1}(z) := \text{ord}_t \partial_z F(t, z),
\label{OR1a}
\end{equation}
and we set
\begin{equation}
\rho_t := \max\{\rho_{t,0}(z), \ \rho_{t,1}(z)\}.
\label{OR1b}
\end{equation}
Then:

(i) The quantity $\rho_t$ is independent of $z$ and is related to the $z$-order $\rho_z$ ($= \rho$) of $F(t, z)$ via the formula
\begin{equation}
\rho_t = \frac{\rho}{2 - \rho}
\label{O18}
\end{equation}
(in particular, if $\rho = 2$, then $\rho_t = \infty$, while $\rho_t$ is finite in the case where $\rho \in [0, 2)$, and $\rho_t = 0$ if and
only if $\rho = 0$).

(ii) We have
\begin{equation}
\rho_{t,0}(z) = \rho_t
\qquad \text{and} \qquad
\rho_{t,1}(z) = \rho_t
\qquad \text{for \;a.e. } z \in \mathbb{C}.
\label{OR2}
\end{equation}

\smallskip

\textit{Proof}. (i) Let us set
\begin{equation}
c_k(z) := \frac{f^{(k)}(z)}{k!},
\qquad\qquad
k \geq 0,
\quad
z \in \mathbb{C}
\label{O4d}
\end{equation}
(so that, in view of \eqref{O4a}, $c_k(0) = c_k$). Then, by substituting \eqref{O4d} in \eqref{O2} we get
\begin{equation}
F(t, z) = \sum_{j \geq 0} \frac{(2j)!}{j!} \, c_{2j}(z) \, t^j
\label{O4e}
\end{equation}
and
\begin{equation}
\partial_z F(t, z) = \sum_{j \geq 0} \frac{(2j)!}{j!} \, c_{2j}'(z) \, t^j = \sum_{j \geq 0} \frac{(2j+1)!}{j!} \, c_{2j+1}(z) \, t^j.
\label{O4f}
\end{equation}
Then, by \eqref{O4e} and \eqref{PREL11} we get
\begin{equation}
\rho_{t,0}(z) = \text{ord}_t F(t, z) = \limsup_j \frac{j \, \ln j}{-\ln|c_{2j}(z) (2j)!/ j!|},
\label{O13}
\end{equation}
which implies
\begin{equation}
\rho_{t,0}(z)
= \limsup_j \frac{1}{\frac{-\ln |c_{2j}(z)|}{j \, \ln j}
- \frac{\ln\left((2j)!\right)}{j \, \ln j} + \frac{\ln\left(j!\right)}{j \, \ln j}}.
\label{O14}
\end{equation}
Now, by Stirling's formula we have that \ $\ln\left((2j)!\right) / (j \, \ln j) \to 2$
and \ $\ln\left(j!\right) / (j \, \ln j) \to 1$ as $j \to \infty$. Thus \eqref{O14} becomes
\begin{equation}
\rho_{t,0}(z)
= \limsup_j \frac{1}{\frac{-\ln |c_{2j}(z)|}{j \, \ln j} - 1}
= \limsup_j \frac{\frac{j \, \ln j}{-\ln|c_{2j}(z)|}}{1 - \frac{j \, \ln j}{-\ln |c_{2j}(z)|}}
= \frac{\rho_0(z)}{2 - \rho_0(z)},
\label{O15}
\end{equation}
where
\begin{equation}
\rho_0(z) := 2 \limsup_j \frac{j \, \ln j}{-\ln |c_{2j}(z)|}.
\label{O15a}
\end{equation}
Likewise, from \eqref{O4f} we have
\begin{equation}
\rho_{t,1}(z) = \text{ord}_t \partial_z F(t, z) = \limsup_j \frac{j \, \ln j}{-\ln |c_{2j+1}(z) (2j+1)!/ j!|}
\label{O16}
\end{equation}
and in the same way as above we obtain
\begin{equation}
\rho_{t,1}(z) = \frac{\rho_1(z)}{2 - \rho_1(z)},
\label{O17}
\end{equation}
where
\begin{equation}
\rho_1(z) : = 2 \limsup_j \frac{j \, \ln j}{-\ln |c_{2j+1}(z)|}.
\label{O17a}
\end{equation}
But, in view of \eqref{O4d}, the quantities $\rho_0(z)$ and $\rho_1(z)$ of \eqref{O15a} and \eqref{O17a} respectively, associated to $f(z)$, are the analogs of the quantities introduced in \eqref{PREL15} and \eqref{PREL16}, associated to $g(z)$. Hence, by Theorem 1 (and Remark 2) we get that
$\rho_0(z) = \rho_1(z) = \rho$ for a.e. $z \in \mathbb{C}$ and, therefore, formulas \eqref{O15} and \eqref{O17}, together with \eqref{OR1b}, imply
\eqref{O18} and \eqref{OR2}.
\hfill $\blacksquare$

\medskip

We will refer to $\rho_t$ of \eqref{OR1b} as the \emph{caloric $t$-order of} $F(t, z)$. Also, let us set
\begin{equation}
\mathcal{E}_0 := \{z \in \mathbb{C} \,:\, \rho_{t,0}(z) = \rho_t\}
\qquad \text{and} \qquad
\mathcal{E}_1 := \{z \in \mathbb{C} \,:\, \rho_{t,1}(z) = \rho_t\},
\label{OR17}
\end{equation}
so that, in view of \eqref{OR1b}, $\mathcal{E}_0 \cup \mathcal{E}_1 = \mathbb{C}$.
Theorem 4 tells us that both $\mathcal{E}_0$ and  $\mathcal{E}_1$ have full measure. Consequently,
\begin{equation}
\mathcal{E}: = \mathcal{E}_0 \cap \mathcal{E}_1
\quad \text{is a subset of $\mathbb{C}$ of full measure}.
\label{OR17a}
\end{equation}

In general it might happen that $\rho_{t,0}(z_0) < \rho_t$ for some $z_0$ and $\rho_{t,1}(z_1) < \rho_t$ for some $z_1 \ne z_0$. For instance, if
$F(t, z) = e^{-\lambda^2 t} \sin(\lambda z)$, then $\rho_t = 1$. However, $F(t, z_0) \equiv 0$ for $z_0 = k\pi/\lambda$, $k \in \mathbb{Z}$, while
$\partial_z F(t, z_1) \equiv 0$ for $z_1 = [k + (1/2)]\pi/\lambda$, $k \in \mathbb{Z}$, and, consequently, $\rho_{t,0}(z_0) = \rho_{t,1}(z_1) = 0$ for those values of $z_0$ and $z_1$.

Next, we present the analog of Theorem 4 regarding the $t$-type of $F(t, z)$.

\medskip

\textbf{Theorem 5.} Let $F(t, z)$ be an entire solution of the heat equation \eqref{P3} whose $z$-order $\rho_z$ ($= \rho$) satisfies
$0 < \rho_z < 2$ (or, equivalently, by \eqref{O18} the caloric $t$-order of $F(t, z)$ satisfies $0 < \rho_t < \infty$). Thinking of $z$ as a
parameter, we consider the types of $F(t, z)$ and $\partial_z F(t, z)$ viewed as functions of $t$, namely
\begin{equation}
\tau_{t,0}(z) := \text{type}_t F(t, z),
\qquad
z \in \mathcal{E}_0
\label{OT1a}
\end{equation}
and
\begin{equation}
\tau_{t,1}(z) := \text{type}_t \partial_z F(t, z),
\qquad
z \in \mathcal{E}_1,
\label{OT1aa}
\end{equation}
where the sets $\mathcal{E}_0$ and $\mathcal{E}_1$ are the full measure sets defined by \eqref{OR17}-\eqref{OR1a}-\eqref{OR1b}. We also set
\begin{equation}
\tau_t := \max\{\tau_{t,0}(z), \ \tau_{t,1}(z)\},
\qquad
z \in \mathcal{E} = \mathcal{E}_0 \cap \mathcal{E}_1.
\label{OT1b}
\end{equation}

Then:

(i) The quantity $\tau_t$ is independent of $z \in \mathcal{E}$ and is related to the $z$-order $\rho_z$ ($= \rho$) and the $z$-type $\tau_z$
($= \tau$) of $F(t, z)$ via the formula
\begin{equation}
\tau_t = \left(1 - \frac{\rho}{2}\right) (2 \rho)^{\frac{\rho}{2-\rho}} \, \tau^{\frac{2}{2-\rho}}.
\label{O24}
\end{equation}
In particular, $\tau_t = 0$ if and only if $\tau = 0$ and $\tau_t = \infty$ if and only if $\tau = \infty$.

(ii) If $\tau_t < \infty$, then
\begin{equation}
\tau_{t,0}(z) = \tau_t
\qquad \text{and} \qquad
\tau_{t,1}(z) = \tau_t
\qquad \text{for \;a.e. } z \in \mathbb{C}.
\label{OT2}
\end{equation}

\smallskip

\textit{Proof}. (i) In view of \eqref{O4d} and \eqref{PREL12} the type $\tau$ of $f(z) = F(0, z)$ is given by
\begin{equation}
\tau = \frac{1}{e \rho} \, \limsup_k k \, |c_k(z)|^{\rho_z/k}
\qquad \text{for any }\; z \in \mathbb{C}.
\label{OT19}
\end{equation}
If we set
\begin{equation}
\tau_0(z) := \frac{2}{e \rho} \, \limsup_j j \, |c_{2j}(z)|^{\rho_z/2j},
\qquad
\tau_1(z) := \frac{2}{e \rho} \, \limsup_j j \, |c_{2j+1}(z)|^{\rho_z/2j},
\label{OT20}
\end{equation}
then, obviously \eqref{OT19} implies
\begin{equation}
\tau = \max\{\tau_0(z), \, \tau_1(z)\}
\qquad \text{for any }\; z \in \mathbb{C},
\label{OT19a}
\end{equation}
while by Theorem 2 we have that, if $\tau_t < \infty$,
\begin{equation}
\tau_0(z) = \tau
\qquad \text{and} \qquad
\tau_1(z) = \tau
\qquad \text{for \;a.e. } z \in \mathbb{C}.
\label{OT20a}
\end{equation}

Next, let us set
\begin{equation}
\hat{\tau}_{t, 0}(z) := \frac{1}{e \rho_t} \, \limsup_j j \, \left|\frac{(2j)! \, c_{2j}(z)}{j!}\right|^{\rho_t/j},
\qquad
z \in \mathbb{C},
\label{O21aa}
\end{equation}
and
\begin{equation}
\hat{\tau}_{t, 1}(z) := \frac{1}{e \rho_t} \, \limsup_j j \, \left|\frac{(2j+1)! \, c_{2j+1}(z)}{j!}\right|^{\rho_t/j},
\qquad
z \in \mathbb{C}.
\label{O21bb}
\end{equation}

Then, formulas \eqref{O4e}, \eqref{O4f}, \eqref{OR17}, \eqref{OT1a}, and \eqref{OT1aa} imply
\begin{equation}
\hat{\tau}_{t, 0}(z) = \tau_{t, 0}(z)
\qquad
\text{for all }\; z \in \mathcal{E}_0,
\label{O21a}
\end{equation}
and
\begin{equation}
\hat{\tau}_{t, 1}(z) = \tau_{t, 1}(z)
\qquad
\text{for all }\; z \in \mathcal{E}_1.
\label{O21b}
\end{equation}
Application of Stirling's asymptotic formula for the factorial to \eqref{O21aa} gives
\begin{align}
\hat{\tau}_{t, 0}(z)
&= \frac{1}{e \rho_t} \, \limsup_j j \, \left(\frac{4j}{e}\right)^{\rho_t} \left(|c_{2j}(z)|^{1/j}\right)^{\rho_t}
\nonumber
\\
&= \frac{4^{\rho_t}}{\rho_t \, e^{\rho_t +1}} \left[ \limsup_j j \left(|c_{2j}(z)|^{1/j}\right)^{\frac{\rho_t}{\rho_t +1}}\right]^{\rho_t +1},
\qquad
z \in \mathbb{C}.
\label{O22a}
\end{align}
Since \eqref{O18} can be written as
\begin{equation}
\frac{\rho_t}{\rho_t +1} = \frac{\rho}{2},
\label{O23}
\end{equation}
in view of \eqref{O18}, \eqref{O23}, and \eqref{OT20}, formula \eqref{O22a} yields
\begin{align}
\hat{\tau}_{t, 0}(z)
&= \frac{4^{\rho_t}}{\rho_t \, e^{\rho_t +1}} \left[ \limsup_j j |c_{2j}(z)|^{\rho/2j}\right]^{\rho_t +1}
\nonumber
\\
&= \frac{4^{\rho_t}}{\rho_t \, e^{\rho_t +1}} \left[\frac{e \rho \tau_0(z)}{2}\right]^{\rho_t +1}
= \left(1 - \frac{\rho}{2}\right) (2 \rho)^{\frac{\rho}{2-\rho}} \, \tau_0(z)^{\frac{2}{2-\rho}},
\qquad
z \in \mathbb{C}.
\label{O30a}
\end{align}
In the same way, starting from \eqref{O21bb} we get
\begin{equation}
\hat{\tau}_{t, 1}(z)
= \left(1 - \frac{\rho}{2}\right) (2 \rho)^{\frac{\rho}{2-\rho}} \, \tau_1(z)^{\frac{2}{2-\rho}},
\qquad
z \in \mathbb{C}.
\label{O30b}
\end{equation}
Therefore \eqref{O24} follows by using \eqref{O21a} and \eqref{O21b} in \eqref{OT1b}, and then invoking \eqref{O30a}, \eqref{O30b} and \eqref{OT19a}.
As for \eqref{OT2}, it follows from \eqref{OT20a}, \eqref{O21a}, \eqref{O21b}, \eqref{O30a}, \eqref{O30b}, and \eqref{OT1b}.
\hfill $\blacksquare$

\medskip

We will refer to $\tau_t$ of \eqref{OT1b} as the \emph{caloric $t$-type of} $F(t, z)$. The question whether \eqref{OT2} remains valid in the case
where $\tau_t = \infty$ remains open.

As an example let us observe that in the special case $\rho_z = \rho = 1$ formula \eqref{O18} implies that $\rho_t = 1$ and then \eqref{O24} yields
$\tau_t = \tau_z^2 = \tau^2$ (e.g., this is the case of the special solution $E_{\lambda}(t, z) = e^{\lambda^2 t + \lambda z}$, where,
clearly, $\rho_z = \rho_t = 1$, $\tau_z = |\lambda|$, and $\tau_t = |\lambda|^2$).

\medskip

\textbf{Remark 7.} From the definition \eqref{OR1b} of the caloric $t$-order $\rho_t$ we know that the $t$-orders
$\rho_{t,0}(z) = \text{ord}_t F(t, z)$ and $\rho_{t,1}(z) = \text{ord}_t \partial_z F(t, z)$ are $\leq \rho_t$ for every
$z \in \mathbb{C}$. However, the $t$-type $\text{type}_t F(t, z)$ of $F(t, z)$ or the $t$-type $\text{type}_t \partial_z F(t, z)$ of
$\partial_z F(t, z)$ can become bigger than the caloric $t$-type $\tau_t$ of $F(t, z)$ for some exceptional values of $z$. For instance, suppose
$f(z) = f_1(z) + f_2(z)$, where the order of $f_1(z)$ is smaller than the order of $f_2(z)$, while the type of $f_1(z)$ is bigger than the type of
$f_2(z)$. Furthermore, let us assume that $f_2(z)$ is an odd function. Then, we can easily construct examples where
$\text{type}_t F(t, 0) > \tau_t$.

\medskip

\textbf{Remark 8.} Let $f(z)$ be an entire function whose exact order is in $[0, 2^-]$. Then,
\begin{equation}
F(t, z) := \int_{-\infty}^{\infty}\frac{1}{2\sqrt{\pi t}} \, e^{-(z - \xi)^2/4t} f(\xi) \, d\xi
\qquad \text{for} \quad \Re(t) > 0
\label{O24a}
\end{equation}
satisfies the heat equation \eqref{P3} for $\Re(t) > 0$. Actually,
\begin{equation}
F(t-t_0, z) = \int_{-\infty}^{\infty}\frac{1}{2\sqrt{\pi (t-t_0)}} \, e^{-(z - \xi)^2/4(t-t_0)} F(t_0, \xi) \, d\xi
\qquad \text{for} \quad \Re(t-t_0) > 0
\label{O24b}
\end{equation}
also satisfies \eqref{P3}.

If $t$ and $z$ are real, with $t > 0$, then we can use the substitution $\eta = (z - \xi)/\sqrt{t}$ or $\eta = (\xi - z)/\sqrt{t}$
in the integral of \eqref{O24a} and obtain
\begin{equation}
F(t, z) = \frac{1}{2\sqrt{\pi}}\int_{-\infty}^{\infty} \, e^{-\eta^2/4} f\left(z + \eta \sqrt{t}\right) \, d\eta
= \frac{1}{2\sqrt{\pi}}\int_{-\infty}^{\infty} \, e^{-\eta^2/4} f\left(z - \eta \sqrt{t}\right) \, d\eta,
\label{O24c}
\end{equation}
but then, due to our assumption for the order and type of $f(z)$, the integral in the right-hand side of \eqref{O24c} is entire in $(t, z)$,
satisfies the heat equation for every $t, z \in \mathbb{C}$ (e.g., by analytic continuation) and it is clear from \eqref{O24c} that $F(0, z) = f(z)$.

\medskip

\textbf{Example 2.} (i) Suppose
\begin{equation}
f(z) = e^{a z^2},
\qquad \text{where }\; a \in \mathbb{C} \setminus \{0\},
\label{O25}
\end{equation}
so that $\rho = 2$ and $\tau = |a| > 0$ (thus the exact order of $f(z)$ is $2$). Here, if $F(t, z)$ satisfies the heat equation and the initial
condition $F(0, z) = f(z)$, then
\begin{equation}
F(t, z) = \frac{1}{\sqrt{1 - 4at}} \, \exp\left(\frac{a z^2}{1 - 4at}\right).
\label{O26}
\end{equation}
Obviously, this $F(t, z)$ is not entire since it has a strong singularity at $t = 1/4a$ (independent of $z$), a combination of an essential
singularity and a square-root branch point (this singularity ``disappears" as $a \to 0$). If $t$ starts at $0$, makes a loop around $1/4a$ and comes
back to $0$, then we get a different value of $F(0, z)$. Thus, if initially $F(0, z) = e^{a z^2}$, then, after one loop we get the other branch of
$F(0, z)$, namely $F(0, z) = -e^{a z^2}$. This counterexample is in agreement with the previous discussion.
Another observation is that for $t \ne 1/4a$ the function $F(t, z)$ is entire in $z$ and its $z$-order is $\rho_z = \rho = 2$, independent of $t$. However, its $z$-type is $\tau_z = |a| |1 - 4at|^{-1}$ (thus depends on $t$). Finally, notice that $F(t, z)$ is never $0$.

(ii) As a variant of the above case consider
\begin{equation}
f(z) = \cos(a z^2) = \frac{e^{ia z^2} + e^{-ia z^2}}{2},
\qquad \text{where }\; a \in \mathbb{C} \setminus \{0\},
\label{O25a}
\end{equation}
so that, again, $\rho = 2$ and $\tau = |a| > 0$. Here, if $F(t, z)$ satisfies the heat equation and the initial condition $F(0, z) = f(z)$, then
\begin{equation}
F(t, z) = \frac{1}{2\sqrt{1 - 4iat}} \exp\left(\frac{ia z^2}{1 - 4iat}\right) + \frac{1}{2\sqrt{1 + 4iat}} \exp\left(\frac{-ia z^2}{1 + 4iat}\right).
\label{O26a}
\end{equation}
Here $F(t, z)$ has a strong singularities at $t = \pm 1/4ia$; furthermore, it has infinitely many zeros. Actually, $z$ is a zero of $F(t, z)$
if and only if
\begin{equation*}
z^2 = (1 + 16 a^2 t^2) \left[\frac{1}{4\pi i a} \ln \left(\frac{1 - 4iat}{1 + 4iat} \right) +\frac{\pi}{2a} \right],
\end{equation*}
where the $\ln(\cdot)$ denotes the multivalued logarithmic function.
\medskip

\textbf{Example 3.} Here we examine the case where the initial condition is
\begin{equation}
f(z) = z^{\alpha},
\qquad
\alpha \in \mathbb{C}.
\label{O27}
\end{equation}
For this $f(z)$ the solution of the heat equation \eqref{P3} cannot be entire in $(t, z)$, unless, of course, $\alpha = m$, a nonnegative integer.

Consider the function
\begin{equation}
F(t, z) := F(t, z; \alpha)
:= \frac{i^{\alpha}}{\Gamma(-\alpha)} \, t^{\alpha/2} \int_0^{\infty} \frac{e^{-\xi^2 + i t^{-1/2} z\xi}}{\xi^{\alpha+1}} \, d\xi,
\qquad
\Re(\alpha) < 0,
\label{O28d}
\end{equation}
where $\Gamma(\cdot)$ is the gamma function. Notice that $F(t, z; \alpha)$ is entire in $z$ for any complex $t \ne 0$ and analytic in $t \ne 0$
for any complex $z$. The
singularity at $t = 0$ it is a combination of a branch point and an essential singularity. Also, it is not hard to
check that $F(t, z)$ satisfies the heat equation for every $z \in \mathbb{C}$, $t \in \mathbb{C}\setminus\{0\}$. Furthermore, if
$t \to 0$ in a way so that $\Im(t^{-1/2} z) \geq 0$, then $F(t, z; \alpha)$ approaches (some branch of) $z^{\alpha}$. However, if
$t \to 0$ in an arbitrary way, then $\lim_{t \to 0} F(t, z; \alpha)$ may not be equal to $z^{\alpha}$. This is the sense in which the initial
condition $F(0, z; \alpha) = z^{\alpha}$ is satisfied.

The integral representation of $F(t, z; \alpha)$ given in \eqref{O28d} makes sense only for $\Re(\alpha) < 0$. We can give other (contour) integral
representations of $F(t, z; \alpha)$ over contours in the complex plane avoiding the positive semiaxis, which are valid for any
$\alpha \in \mathbb{C}$. But, instead of doing that, we continue the analysis of $F(t, z; \alpha)$ as follows.

We observe that $F(t, z) = F(t, z; \alpha)$ of \eqref{O28d} can be written as
\begin{equation}
F(t, z) = \frac{i^{\alpha}  t^{\alpha/2}}{\Gamma(-\alpha)} \, h\left(\frac{it^{-1/2}z}{2}\right),
\quad \text{where}\quad
h(x) := h(x; \alpha) := \int_0^{\infty} \frac{e^{-\xi^2 + 2x\xi}}{\xi^{\alpha+1}} \, d\xi,
\label{O28c}
\end{equation}
where the function $h(x)$ of \eqref{O28c} satisfies the \emph{Hermite equation}, namely
\begin{equation}
u''(x) - 2x u'(x) = -2\alpha u(x),
\label{O29}
\end{equation}
with
\begin{equation}
h(0) = \frac{1}{2} \, \Gamma\left(-\frac{\alpha}{2}\right)
\qquad \text{and} \qquad
h'(0) = \Gamma\left(\frac{1-\alpha}{2}\right).
\label{O29a}
\end{equation}
Clearly, every solution of \eqref{O29} is entire in $x$. Actually, the general solution of \eqref{O29} can be expressed as
\begin{equation}
u(x) = \sum_{n=0}^{\infty} a_n x^k,
\label{O29b}
\end{equation}
where the coefficients $a_k$ satisfy the recursion
\begin{equation}
\frac{a_{n+2}}{a_n} = \frac{2(n - \alpha)}{(n+1)(n+2)}.
\end{equation}
We can single out two linearly independent solutions of \eqref{O29} in a convenient way. By taking $a_0 = 1$ and $a_1 = 0$ we obtain the solution
$u_e(x) = u_e(x; \alpha)$ which is even in $x$ and satisfies $u_e(0) = 1$, while by taking $a_0 = 0$ and $a_1 = 1$ we obtain the solution
$u_o(x) = u_o(x; \alpha)$ which is odd in $x$ and satisfies $u_o'(0) = 1$. Thus,
\begin{equation}
u_e(x) = 1 + \sum_{k=1}^{\infty} (-1)^k\frac{2^k \alpha (\alpha - 2) \cdots \big(\alpha - 2(k-1)\big)}{(2k)!} \, x^{2k}
\label{O29f}
\end{equation}
and
\begin{equation}
u_o(x) = x + \sum_{k=1}^{\infty} (-1)^k\frac{2^k (\alpha - 1)(\alpha - 3) \cdots (\alpha - 2k + 1)}{(2k+1)!} \, x^{2k+1}.
\label{O29g}
\end{equation}
Evidently (due to the analytic dependence on the parameter $\alpha$), $u_e(x; \alpha)$ and $u_o(x; \alpha)$ are also entire in $\alpha$. Furthermore,
in view of \eqref{O29a},
\begin{equation}
h(x) = h(x; \alpha) = \frac{1}{2} \, \Gamma\left(-\frac{\alpha}{2}\right) u_e(x; \alpha) + \Gamma\left(\frac{1-\alpha}{2}\right) u_o(x; \alpha).
\label{O28e}
\end{equation}
An additional consequence of formula \eqref{O28e} is that, by analytic continuation in $\alpha$, we get the meromorphic extension of $h(x; \alpha)$,
which we also denote by $h(x; \alpha)$, for all complex $\alpha$.

Using \eqref{O28e} in \eqref{O28c} yields
\begin{equation}
F(t, z; \alpha) = \frac{i^{\alpha}  t^{\alpha/2}}{\Gamma(-\alpha)}\left[
\frac{1}{2} \, \Gamma\left(-\frac{\alpha}{2}\right) u_e\left(\frac{it^{-1/2}z}{2}; \alpha\right)
+ \Gamma\left(\frac{1-\alpha}{2}\right) u_o\left(\frac{it^{-1/2}z}{2}; \alpha\right)
\right].
\label{O30}
\end{equation}
With the help of the well-known Legendre's duplication formula for the Gamma function, namely the formula
\begin{equation*}
\sqrt{\pi} \, \Gamma(2z) = 2^{2z-1} \Gamma(z) \Gamma\left(z + \frac{1}{2}\right),
\end{equation*}
equation \eqref{O30} simplifies as
\begin{equation}
F(t, z; \alpha) = \sqrt{\pi} (2i)^{\alpha} t^{\alpha/2}\left[
\frac{1}{\Gamma\left(\frac{1-\alpha}{2}\right)} \,  u_e\left(\frac{it^{-1/2}z}{2}; \alpha\right)
+ \frac{2}{\Gamma(-\alpha/2)} u_o\left(\frac{it^{-1/2}z}{2}; \alpha\right)
\right].
\label{O31}
\end{equation}
It follows that for $t \ne 0$ the function $F(t, z; \alpha)$ is entire in $z$ and $\alpha$, while for any $z$ and $\alpha$ it is analytic in $t$, except for $t=0$, where it may have a strong singularity (branch point combined with an essential singularity). Thus, if $t$ makes a loop
around $0$ in the complex $t$-plane, then we may arrive at a different value of $F(t, z; \alpha)$. Notice, however, that if $m \geq 0$ is an integer, $F(t, z; m)$ becomes the $m$-th caloric polynomial, i.e.
\begin{equation}
F(t, z; m) = P_m(t,z).
\label{O32}
\end{equation}
Generically, the $z$-order of $F(t, z; \alpha)$ is $2$, while its $z$-type depends on $t$.

The case $\alpha = -1$ is of particular interest. The solutions $u_e(x; \alpha)$ and $u_o(x; \alpha)$ become respectively
\begin{equation}
u_e(x; -1) = e^{x^2}
\qquad \text{and} \qquad
u_o(x; -1) = e^{x^2} \int_0^x e^{-\xi^2} d\xi
\label{O33}
\end{equation}
for every $x \in \mathbb{C}$. Then, with the help of \eqref{O33} and \eqref{O31} we can construct the solution of the heat equation
\begin{equation}
F(t, z) = \frac{i}{\sqrt{t}} \, e^{-z^2/4t} \int_{-\infty}^{-iz/2\sqrt{t}} e^{-\zeta^2} d\zeta
= \frac{i}{\sqrt{t}} \, e^{-z^2/4t} \left(\frac{\sqrt{\pi}}{2} + \int_0^{-iz/2\sqrt{t}} e^{-\zeta^2} d\zeta\right),
\label{O34}
\end{equation}
where the first contour integral is taken over a contour which approaches the negative real axis at $-\infty$.
Formula \eqref{O34} implies that
\begin{equation}
F(t, z) \to \frac{1}{z}
\quad \text{as }\; t \to 0
\quad
\text{in certain directions}.
\label{O35}
\end{equation}

From $F(t, z)$ of \eqref{O34} we can also obtain the solution of the heat equation with
``initial condition" $f(z) = \ln z$ as $\int_1^z F(t, \zeta) d\zeta$ (the quotations here remind us that the initial condition is satisfied in a
certain sense). Furthermore, the solutions of the heat equation with ``initial conditions" $f(z) = z^{-m}$, $m \in \mathbb{N}$,
can be also obtained from $F(t, z)$ by differentiating it $m-1$ times with respect to $z$ (or $(m-1)/2$ times with respect to $t$, in $m$ is odd).

\section{The zeros of entire caloric functions}
We start with a result stating that the multiple zeros of a (nontrivial) entire caloric function $F(t, z)$, viewed as a function of $z$, cannot
accumulate in $\mathbb{C}^2$.

\medskip

\textbf{Theorem 6.} Suppose $F(t, z) \not\equiv 0$ is entire in $(t, z)$ and satisfies the heat equation \eqref{P3}. If
\begin{equation}
F(t^{\ast}, z^{\ast}) = \partial_z F(t^{\ast}, z^{\ast}) = 0
\qquad \text{for some} \quad
(t^{\ast}, z^{\ast}) \in \mathbb{C}^2,
\label{P14}
\end{equation}
then there is a ($\mathbb{C}^2$-open) neighborhood $U$ of $(t^{\ast}, z^{\ast})$ such that
\begin{equation}
|F(t, z)| + |\partial_z F(t, z)| > 0
\qquad \text{for every} \quad
(t, z) \in U \setminus\{(t^{\ast}, z^{\ast})\}.
\label{P15}
\end{equation}

\smallskip

\textit{Proof}. Without loss of generality and for typographical convenience we take $t^{\ast} = z^{\ast} = 0$.

If the statement of the theorem is false, then there exists a sequence of points $(t_n, z_n) \ne (0, 0)$, $n = 1, 2, \ldots$, such that
$(t_n, z_n) \to (0,0)$
and
\begin{equation}
F(t_n, z_n) = \partial_z F(t_n, z_n) = 0
\qquad \text{for every} \quad
n \geq 1.
\label{P18}
\end{equation}

Let $\mu$ be the smallest value of $m$ such that $\partial_z^m F(0, 0) \ne 0$ (since $F(t, z) \not\equiv 0$, formula \eqref{P6} guarantees that $\mu$
exists; of course, due to the assumption \eqref{P14} we have $\mu \geq 2$). Then, by taking $t_0 = z_0 = 0$ in \eqref{P6} we get
\begin{equation}
F(t, z) = \frac{\partial_z^{\mu} F(0, 0)}{\mu!} \, P_{\mu}(t, z) \,
+ \sum_{m=\mu+1}^{\infty} \frac{\partial_z^m F(0, 0)}{m!} \, P_m(t, z),
\label{P16c}
\end{equation}
which also implies
\begin{equation}
\partial_z F(t, z) = \frac{\partial_z^{\mu} F(0, 0)}{\mu!} \, \partial_z P_{\mu}(t, z) \,
+ \sum_{m=\mu+1}^{\infty} \frac{\partial_z^m F(0, 0)}{m!} \, \partial_z P_m(t, z).
\label{P16d}
\end{equation}
In view of \eqref{P7}, formulas \eqref{P16c} and \eqref{P16d} imply that, given an open ball
$B \subset \mathbb{C}^2$ centered at $(0, 0)$ there is an $C > 0$ (i.e. depending only on $B$) such that for every $(t, z) \in B$ we have
\begin{equation}
\left|F(t, z) - \frac{\partial_z^{\mu} F(0, 0)}{\mu!} \, P_{\mu}(t, z)\right|
\leq C \max\left\{|z|^{\mu+1}, |z|^{\mu-1} |t|, \ldots, |t|^{\lfloor(\mu+2)/2\rfloor}\right\}
\label{P16e}
\end{equation}
and
\begin{equation}
\left|\partial_z F(t, z) - \frac{\partial_z^{\mu} F(0, 0)}{\mu!} \, \partial_z P_{\mu}(t, z)\right|
\leq C \max\left\{|z|^{\mu}, |z|^{\mu-2} |t|, \ldots, |t|^{\lfloor(\mu+1)/2\rfloor}\right\}.
\label{P16f}
\end{equation}
Now, given $B$ there is an $n_0$ such that $(t_n, z_n) \in B$ for all $n \geq n_0$. Hence, by using \eqref{P18} in \eqref{P16e} and \eqref{P16f}
we get
\begin{equation}
\left|P_{\mu}(t_n, z_n)\right| \leq C' \max\left\{|z_n|^{\mu+1}, |z_n|^{\mu-1} |t_n|, \ldots, |t_n|^{\lfloor(\mu+2)/2\rfloor}\right\}
\quad \text{for every }
n \geq n_0
\label{P16g}
\end{equation}
and
\begin{equation}
\left|\partial_z P_{\mu}(t_n, z_n)\right| \leq C' \max\left\{|z_n|^{\mu}, |z_n|^{\mu-2} |t_n|, \ldots, |t_n|^{\lfloor(\mu+1)/2\rfloor}\right\}
\quad \text{for every }
n \geq n_0,
\label{P16h}
\end{equation}
where for typographical convenience we have set
\begin{equation}
C' :=  \frac{\mu! C}{|\partial_z^{\mu} F(0, 0)|}.
\label{P16i}
\end{equation}

Let us consider the case where $\mu = 2l$. Then, substituting \eqref{P12} in \eqref{P16g} yields
\begin{equation}
|z_n^2 + \rho_{\mu,1} \, t_n| \cdots |z_n^2 + \rho_{\mu,l} \, t_n|
\leq C' \max\left\{|z_n|^{2l+1}, |z_n|^{2l-1} |t_n|, \ldots, |z_n| \, |t_n|^l, |t_n|^{l+1}\right\}
\label{P17a}
\end{equation}
for every $n \geq n_0$, while, in view of \eqref{P7d}, substituting \eqref{P13} in \eqref{P16h} yields
\begin{equation}
|z_n||z_n^2 + \rho_{\mu-1,1} \, t_n| \cdots |z_n^2 + \rho_{\mu-1,l-1} \, t_n|
\leq \frac{C'}{\mu} \max\left\{|z_n|^{2l}, |z_n|^{2l-2} |t_n|, \ldots, |z_n|^2 |t_n|^{l-1}, |t_n|^l\right\}
\label{P17b}
\end{equation}
for every $n \geq n_0$.

If $z_n = 0$ (hence $t_n \ne 0$), then \eqref{P17a} becomes $|\rho_{\mu,1} \cdots \rho_{\mu,l} \, t_n^l| \leq C' |t_n|^{l+1}$, which, in view of
\eqref{P12a} and the assumption $t_n \to 0$, cannot be satisfied for any sufficiently large $n$. Thus, without loss of generality we can assume
$z_n \ne 0$.

We set
\begin{equation}
\lambda_n := \frac{t_n}{z_n^2}.
\label{P19}
\end{equation}
Then, \eqref{P17a} and \eqref{P17b} become respectively
\begin{equation}
|1 + \rho_{\mu,1} \, \lambda_n| \cdots |1+ \rho_{\mu,l} \, \lambda_n|
\leq C' |z_n| \max\left\{1, |\lambda_n|, \ldots, |\lambda_n|^l, |\lambda_n|^{l+1} |z_n|\right\}
\label{P19a}
\end{equation}
and
\begin{equation}
|1 + \rho_{\mu-1,1} \, \lambda_n| \cdots |1 + \rho_{\mu-1,l-1} \, \lambda_n|
\leq \frac{C'}{\mu} |z_n| \max\left\{1, |\lambda_n|, \ldots, |\lambda_n|^l\right\}
\label{P19b}
\end{equation}
for every sufficiently large $n$.

If $|\lambda_n|$ becomes arbitrarily large, then, in view of \eqref{P12a}, formula \eqref{P19a} should imply that there is a constant $C'' >0$
such that
\begin{equation}
|\lambda_n|^l \leq C'' |z_n| \max\left\{|\lambda_n|^l, |\lambda_n|^{l+1} |z_n|\right\}
= C'' \max\left\{|\lambda_n|^l |z_n|, |\lambda_n|^l |t_n|\right\},
\label{P19c}
\end{equation}
i.e. $1 \leq C'' \max\left\{|z_n|, |t_n|\right\}$, which is, obviously, impossible since $z_n, t_n \to 0$. Therefore, the sequence $\lambda_n$ must be bounded and, hence, by \eqref{P19a} and \eqref{P19b} there must exist a constant $M > 0$ such that
\begin{equation}
|1 + \rho_{\mu,1} \, \lambda_n| \cdots |1+ \rho_{\mu,l} \, \lambda_n| \leq M |z_n|
\label{P19d}
\end{equation}
and
\begin{equation}
|1 + \rho_{\mu-1,1} \, \lambda_n| \cdots |1 + \rho_{\mu-1,l-1} \, \lambda_n| \leq M |z_n|
\label{P19f}
\end{equation}
for every sufficiently large $n$. Let $\lambda_{n_k}$ be a convergent subsequence of the sequence $\lambda_n$, with
$\lim \lambda_{n_k} = \lambda \in \mathbb{C}$. However, if take limits in \eqref{P19d} and \eqref{P19f} as $n_k \to \infty$ we get
\begin{equation}
|1 + \rho_{\mu,1} \, \lambda| \cdots |1+ \rho_{\mu,l} \, \lambda| = 0
\label{P19g}
\end{equation}
and
\begin{equation}
|1 + \rho_{\mu-1,1} \, \lambda| \cdots |1 + \rho_{\mu-1,l-1} \, \lambda| = 0
\label{P19h}
\end{equation}
which contradict \eqref{P12a}.

The remaining case is $\mu = 2l+1$. Here, by substituting \eqref{P13} in \eqref{P16g} we get
\begin{equation}
|z_n||z_n^2 + \rho_{\mu,1} \, t_n| \cdots |z_n^2 + \rho_{\mu,l} \, t_n|
\leq C' \max\left\{|z_n|^{2l+2}, |z_n|^{2l} |t_n|, \ldots, \, |t_n|^{l+1}\right\}
\label{P20a}
\end{equation}
for every $n \geq n_0$, while, in view of \eqref{P7d}, substituting \eqref{P12} in \eqref{P16h} yields
\begin{equation}
|z_n^2 + \rho_{\mu-1,1} \, t_n| \cdots |z_n^2 + \rho_{\mu-1,l} \, t_n|
\leq \frac{C'}{\mu} \max\left\{|z_n|^{2l+1}, |z_n|^{2l-1} |t_n|, \ldots, |z_n| \, |t_n|^l, |t_n|^{l+1}\right\}
\label{P20b}
\end{equation}
for every $n \geq n_0$.

By proceeding in the same manner as in the case $\mu = 2l$, we again arrive at a contradiction. Therefore, our assumption of the
existence of the sequence $(t_n, z_n)$ is false.
\hfill $\blacksquare$

\medskip

An immediate consequence of Theorem 6 is that, if $F(t, z) \not\equiv 0$ is entire in $(t, z)$ and satisfies the heat equation, then the set
\begin{equation}
\mathcal{M}_F := \{(t, z) \in \mathbb{C}^2 \, : \, F(t, z) = \partial_z F(t, z) = 0\}
\label{P21}
\end{equation}
is discrete in $\mathbb{C}^2$. For example, if $F$ is a caloric polynomial, then by \eqref{P12} and \eqref{P13} we have
\begin{equation}
\mathcal{M}_{P_m} = \{(0, 0)\},
\quad
m \geq 2,
\qquad \text{while} \qquad
\mathcal{M}_{P_0} = \mathcal{M}_{P_1} = \emptyset.
\label{P21a}
\end{equation}

One peculiar consequence of Theorem 6 (together with the fact \cite{K} that zeros of entire functions of two or more complex variables are
never isolated) is that if an entire function $A(t, z)$ can be written as
\begin{equation}
A(t, z) = A_1(t, z)^2 A_2(t, z),
\label{P23}
\end{equation}
where $A_1(t, z)$ and $A_2(t, z)$ are entire and $A_1(t_0, z_0) = 0$ for some point $(t_0, z_0) \in \mathbb{C}^2$, then $A(t, z)$ cannot
satisfy the heat equation \eqref{P3}.

Finally, let us notice that the analog to Theorem 6 in the case where $F(t, z)$ is viewed as a function of $t$ does not hold. For example if we consider the caloric function $F(t, z) = e^{-\lambda^2 t} \sin(\lambda z)$, then $F(t, 0) \equiv 0$, hence $\partial_t^j F(t, 0) \equiv 0$ for every
$j \in \mathbb{N}$.

\subsection{The dynamics of the zeros}
The following lemma is well known (see, e.g., \cite{C-S-V} Lemma 2.3). Since its proof is very short, we include it here for the sake of
completeness.

\medskip

\textbf{Lemma 3.} Let $g(z)$ be analytic in a domain $D$ of $\mathbb{C}$. If $z_0 \in D$ is such that $g(z_0) \ne 0$ and we set
\begin{equation*}
G(z) := (z - z_0) g(z),
\end{equation*}
then
\begin{equation*}
\frac{G''(z_0)}{G'(z_0)} = 2 \, \frac{g'(z_0)}{g(z_0)}.
\end{equation*}

\smallskip

\textit{Proof}. For $z \in D$ We have
\begin{equation*}
\frac{G''(z)}{G'(z)} = \frac{(z - z_0) g''(z) + 2g'(z)}{(z - z_0) g'(z) + g(z)}
\end{equation*}
and the statement follows by setting $z = z_0$.
\hfill $\blacksquare$

\medskip

\textbf{Corollary 2.} Suppose $G(z)$ is an entire function with zeros $z_0, z_1, z_2, \ldots\,$, where $z_0$ is a simple zero of $G(z)$. Furthermore,
let us also asssume (essentially without loss of generality) that $G(0) \ne 0$.

(i) If $\Sigma_{k \geq 0} |z_k|^{-1} < \infty$ and
\begin{equation}
G(z) = C e^{A z} \prod_{k \geq 0} \left(1 - \frac{z}{z_k}\right),
\label{LZ1a}
\end{equation}
where $A$ and $C \ne 0$ are complex constants, then
\begin{equation}
\frac{G''(z_0)}{G'(z_0)} = 2A + 2\sum_{k \geq 1} \frac{1}{z_0 - z_k}.
\label{LZ2a}
\end{equation}

(ii) If $\Sigma_{k \geq 0} |z_k|^{-2} < \infty$ and
\begin{equation}
G(z) = C e^{A z} \prod_{k \geq 0} \left(1 - \frac{z}{z_k}\right) e^{z/z_k},
\label{LZ1b}
\end{equation}
where, again, $A$ and $C \ne 0$ are complex constants, then
\begin{equation}
\frac{G''(z_0)}{G'(z_0)} = 2A + 2\sum_{k \geq 1} \left(\frac{1}{z_0 - z_k} + \frac{1}{z_k}\right).
\label{LZ2b}
\end{equation}

\smallskip

\textit{Proof}. Set
\begin{equation}
g(z) := \frac{G(z)}{z - z_0}
\label{LZ3}
\end{equation}
(thus $G(0) \ne 0$ implies $g(0) \ne 0$). Then, by Lemma 3 we have
\begin{equation}
\frac{G''(z_0)}{G'(z_0)} = 2 \, \frac{g'(z_0)}{g(z_0)}.
\label{LZ4}
\end{equation}
Using \eqref{LZ1a} in \eqref{LZ3} yields
\begin{equation}
g(z) = c e^{A z} \prod_{k \geq 1} \left(1 - \frac{z}{z_k}\right),
\label{LZ5}
\end{equation}
where $c \ne 0$. Thus, the Mittag-Leffler expansion of $g'(z)/g(z)$ is \cite{H}
\begin{equation*}
\frac{g'(z)}{g(z)} = A + \sum_{k \geq 1} \frac{1}{z - z_k}.
\end{equation*}
Therefore, formula \eqref{LZ2a} follows by setting $z = z_0$ in the above formula and substituting in \eqref{LZ4}.

The proof of formula \eqref{LZ2b} is very similar.
\hfill $\blacksquare$

\medskip

Notice that formula \eqref{LZ2b} differs from \eqref{LZ2a} only if $\Sigma_{k \geq 0} |z_k|^{-1} = \infty$.

Let us now consider the set
\begin{equation}
\Gamma := \left\{F(t, z) = 0\right\} := \left\{(t, z) \in \mathbb{C}^2 \, : \, F(t, z) = 0\right\} .
\label{Z00}
\end{equation}
If $\Gamma$ is empty, that is if $F(t, z)$ is never $0$, then by the Hadamard Factorization Theorem and by Theorems 3 and 4 of the previous section
it follows that (recall \eqref{P9})
\begin{equation}
F(t, z) = c E_{\lambda}(t, z) = c \, e^{\lambda^2 t + \lambda z},
\label{Z1}
\end{equation}
for some constants $c, \lambda \in \mathbb{C}$, with $c \ne 0$.

Suppose now that $F(t, z)$ is a nontrivial entire solution of the heat equation, which is not of the form \eqref{Z1} and, hence, it vanishes for some values of $t$ and $z$. Then, $\Gamma$ is a nonempty set in $\mathbb{C}^2$ (and, of course, $\Gamma \ne \mathbb{C}^2$). If $F(t, z)$ is \emph{irreducible}, namely it cannot be written as
\begin{equation}
F(t, z) = A_1(t, z) \, A_2(t, z),
\label{Z2}
\end{equation}
where both $A_1(t, z)$ and $A_2(t, z)$ are entire and assume the value $0$, then we can say that $\Gamma$ is a ``curve" in $\mathbb{C}^2$. Otherwise
$\Gamma$ is a union of such component-curves (as we have seen, a consequence of Theorem 6 is that multiple components do not exist). For example, if
$F$ is a caloric polynomial, then, in view of \eqref{P12} and \eqref{P13} we have
\begin{equation}
\Gamma = \bigcup_{j=1}^{l} \{z^2 + \rho_{m,j} \, t = 0\},
\qquad \text{if} \quad
m = 2l
\label{Z3a}
\end{equation}
and
\begin{equation}
\Gamma = \{z = 0\} \cup \bigcup_{j=1}^{l} \{z^2 + \rho_{m,j} \, t = 0\},
\qquad \text{if} \quad
m = 2l+1.
\label{Z3b}
\end{equation}

Let $z_1(t), z_2(t), \ldots$  be the zeros of $F(t, z)$. These zeros can be seen as branches of a global analytic function, say $Z(T)$ defined on a Riemann surface which can be identified with $\Gamma$. In general, this will be an infinitely sheeted surface. The ramification
points of $\Gamma$ are the points $(t^{\ast}, z^{\ast})$ satisfying \eqref{P14} and Theorem 6 assures us that they form a discrete set in $\mathbb{C}^2$. Thus, for every zero $z_k(t)$ we have that $\partial_z F(t, z_k(t)) \ne 0$ for a.a. $t \in \mathbb{C}$, where here ``a.a." means
``almost all," namely all except for a discrete subset of $\mathbb{C}$. Then, by differentiating (implicitly) $F(t, z_k(t)) = 0$ with respect to $t$ we get
\begin{equation}
\partial_t F(t, z_k(t)) + \partial_z F(t, z_k(t)) z_k'(t) = 0
\label{Z4}
\end{equation}
or, in view of \eqref{P3},
\begin{equation}
z_k'(t) = -\frac{\partial_t F(t, z_k(t))}{\partial_z F(t, z_k(t))} = -\frac{\partial_z^2 F(t, z_k(t))}{\partial_z F(t, z_k(t))}
\qquad\text{for a.a. } t \in \mathbb{C}.
\label{Z5}
\end{equation}
Thus, if for $F(t, z)$ we have that $\rho_z < 1$, then we can apply Corollary 2 to \eqref{Z5} and get
\begin{equation}
z'_k(t) = -2 \sum_{j \ne k} \frac{1}{z_k(t) - z_j(t)}
\qquad\text{for a.a. } t \in \mathbb{C}.
\label{Z6}
\end{equation}
The above derivation of \eqref{Z6} is an imitation of the derivation of \eqref{RT5} as presented in \cite{C-S-V}.

As an application, let us consider the system of ordinary differential equations
\begin{equation}
z'_k(t) = -2 \sum_{j \ne k} \frac{1}{z_k(t) - z_j(t)},
\qquad
1 \leq k \leq N,
\label{SZ6a}
\end{equation}
with initial condition
\begin{equation}
z_k(0) = a_k
\qquad
1 \leq k \leq N,
\label{SZ6b}
\end{equation}
where $a_1, \ldots, a_N$ are distinct non-zero complex numbers. To solve this system, we form the polynomial
\begin{equation}
f(z) := \left(1 - \frac{z}{a_1}\right) \cdots \left(1 - \frac{z}{a_N}\right) = 1 + A_1z + \cdots + A_N z^N.
\label{SZ7}
\end{equation}
Then, the solution $z_1(t), \ldots, z_N(t)$ of the system is the set of zeros of the polynomial in $(t, z)$ given by
\begin{equation}
F(t, z) := 1 + \sum_{k=1}^N A_k P_k(t, z),
\label{SZ8}
\end{equation}
where $P_k(t, z)$ is the $k$-th caloric polynomial.

We expect that this application extends to the infinite case (i.e. $N = \infty$) under the restriction that $\sum_k |a_k|^{-1} < \infty$.

\subsubsection{Even caloric functions}
Suppose that the initial condition $f(z)$ satisfies $f(-z) = f(z)$, i.e. is even. Then, by formula \eqref{O2} we have that the solution $F(t, z)$ of
the heat equation also satisfies $F(t, -z) = F(t, z)$. Hence $F(t, z) = \Phi(t, z^2)$, where $\Phi(t, \mu)$ is entire in $(t, \mu)$ and if $\rho_z$ is the $z$-order of $F(t, z)$, then the order of $\Phi(t, \mu)$ with respect to $\mu$ (the $\mu$-order) is $\rho_z/2$. Furthermore, $\Phi(t, \mu)$
satisfies the heat-type equation
\begin{equation}
\partial_t \Phi(t, z) = 4\mu \partial_{\mu}^2 \Phi(t, \mu) + 2\partial_{\mu} \Phi(t, \mu).
\label{Z7}
\end{equation}
Let $\pm z_1(t), \pm z_2(t), \ldots$  be the zeros of $F(t, z)$. Then the zeros of $\Phi(t, \mu)$ are
$\mu_1(t) = z_1(t)^2, \mu_2(t) = z_2(t)^2, \ldots$, and by imitating the derivation of \eqref{Z5} we now have, in view of \eqref{Z7}
\begin{equation}
\mu_k'(t) = -\frac{\partial_t \Phi(t, \mu_k(t))}{\partial_\mu \Phi(t, \mu_k(t))}
= -4 \mu_k(t) \frac{ \partial_{\mu}^2 \Phi(t, \mu_k(t))}{\partial_\mu \Phi(t, \mu_k(t))} - 2
\qquad\text{for a.a. } t \in \mathbb{C}.
\label{Z8}
\end{equation}
If $\rho_z < 2$, then the $\mu$-order of $\Phi(t, \mu)$ is less than $1$. We can, therefore apply Corollary 2 to \eqref{Z8} and get
\begin{equation}
\mu'_k(t) = -2 - 8 \mu_k(t) \sum_{j \ne k} \frac{1}{\mu_k(t) - \mu_j(t)}.
\label{Z9}
\end{equation}

\subsubsection{Odd caloric functions}
Now, suppose that the initial condition $f(z)$ satisfies $f(-z) = -f(z)$, i.e. is odd. Then, by formula \eqref{O2} we have that the solution
$F(t, z)$ of the heat equation also satisfies $F(t, -z) = -F(t, z)$. Hence $F(t, z) = z\Psi(t, z^2)$, where $\Psi(t, \mu)$ is entire in $(t, \mu)$ and if $\rho_z$ is the $z$-order of $F(t, z)$, then the $\mu$-order of $\Psi(t, \mu)$ is $\rho_z/2$. Furthermore, $\Psi(t, \mu)$
satisfies the heat-type equation
\begin{equation}
\partial_t \Psi(t, z) = 4\mu \partial_{\mu}^2 \Psi(t, \mu) + 6\partial_{\mu} \Psi(t, \mu).
\label{ZO7}
\end{equation}
Let $z_0(t) \equiv 0, \pm z_1(t), \pm z_2(t), \ldots$  be the zeros of $F(t, z)$. Then the zeros of $\Psi(t, \mu)$ are
$\mu_1(t) = z_1(t)^2, \mu_2(t) = z_2(t)^2, \ldots$, and by imitating the derivation of \eqref{Z5} we now have, in view of \eqref{ZO7}
\begin{equation}
\mu_k'(t) = -\frac{\partial_t \Psi(t, \mu_k(t))}{\partial_\mu \Psi(t, \mu_k(t))}
= -4 \mu_k(t) \frac{ \partial_{\mu}^2 \Psi(t, \mu_k(t))}{\partial_\mu \Psi(t, \mu_k(t))} - 6
\qquad\text{for a.a. } t \in \mathbb{C}.
\label{ZO8}
\end{equation}
If $\rho_z < 2$, then the $\mu$-order of $\Psi(t, \mu)$ is less than $1$. We can, therefore apply Corollary 1 to \eqref{ZO8} and get
\begin{equation}
\mu'_k(t) = -6 - 8 \mu_k(t) \sum_{j \ne k} \frac{1}{\mu_k(t) - \mu_j(t)}.
\label{ZO9}
\end{equation}

\subsubsection{Some general examples}
Let us start with an observation.

\medskip

\textbf{Observation 1.} Suppose $F(t, z)$ satisfies the heat equation \eqref{P3}, with $F(0, z) = f(z)$. Then
\begin{equation}
G(t, z) := F(t, z + 2 \lambda t)\, E_{\lambda}(t, z) = F(t, z + 2 \lambda t)\, e^{\lambda^2 t + \lambda z},
\label{L10}
\end{equation}
where $\lambda$ is a complex constant, also satisfies \eqref{P3}, with $G(0, z) = e^{\lambda z} f(z)$.

\medskip

It is quite easy to check the validity of Observation 1. As an example, let us take $f(z) = z^m$, where $m$ is a positive integer. Then
$F(t, z) = P_m(t, z)$, the $m$-th caloric polynomial, and
\begin{equation}
G(t, z) = e^{\lambda^2 t + \lambda z} P_m(t, z + 2 \lambda t),
\label{L11}
\end{equation}
is the solution of the heat equation with initial condition $G(0, z) = z^m e^{\lambda z}$.

(i) Let $F(t, z)$ and $G(t, z)$ be entire solutions of the heat equation \eqref{P3} with initial conditions
\begin{equation}
f(z) = \prod_k \left(1 - \frac{z}{a_k}\right)
\qquad \text{and} \qquad
g(z) = e^{\lambda z} \prod_k \left(1 - \frac{z}{a_k}\right)
\label{Z10a}
\end{equation}
respectively, where the order of the product $\Pi_k [1 - (z / a_k)]$ is $\sigma < 1$ (in other words, there is an $\alpha < 1$ such that
$\Sigma_k |a_k|^{-\alpha} < \infty$). Then, as we have seen the $z$-order of $F(t, z)$ is $\sigma$. It follows that if
$z_1(t), z_2(t), \ldots$  are the zeros of $F(t, z)$, then $\Sigma_k' |z_k(t)|^{-\alpha} < \infty$ for some $\alpha < 1$ (the prime on the sum
indicates that we omit the $z_k(t)$'s which are equal to $0$).

The relation of $G(t, z)$ and $F(t, z)$ is given by \eqref{L10}. Thus, if $w_1(t), w_2(t), \ldots$  are the zeros of $G(t, z)$ (viewed as a function of $z$), then \eqref{L10} implies that
\begin{equation}
w_k(t) = z_k(t) - 2\lambda t,
\qquad
k \geq 1,
\label{Z11}
\end{equation}
and it follows that $\Sigma_k' |w_k(t)|^{-\alpha} < \infty$.

Finally, since \eqref{L11} implies $w_k'(t) = z_k'(t) - 2\lambda$, while $z_k(t)$, $k \geq 1$, satisfies \eqref{Z6}, we have
\begin{equation}
w'_k(t) = -2\lambda -2 \sum_{j \ne k} \frac{1}{w_k(t) - w_j(t)}
\qquad\text{for a.a. } t \in \mathbb{C}.
\label{Z12}
\end{equation}

\medskip

(ii) Let $F(t, z)$ be an entire solution of the heat equation \eqref{P3} with initial condition
\begin{equation}
f(z) = e^{\lambda z} z^d\prod_k \left(1 - \frac{z}{a_k}\right) e^{z/a_k},
\label{Z13}
\end{equation}
where, $d \geq 0$ is an integer and the order $\rho$ of $f(z)$ is $< 2$. Then, by applying Corollary 2(ii) we can get that the zeros
$z_1(t), z_2(t), \ldots$  of $F(t, z)$ satisfy
\begin{equation*}
z'_k(t) = -2\lambda -2 \sum_{j \ne k} \left[\frac{1}{z_k(t) - z_j(t)} + \frac{1}{z_j(t)}\right]
\qquad\text{for a.a. } t \in \mathbb{C}.
\end{equation*}

%
%
%
%
%

\end{document}